\documentclass[
final
]{dmtcs-episciences}
\usepackage[utf8]{inputenc} 
\usepackage{amsthm} 
\usepackage{amsmath,amsfonts}
\usepackage{bbold}
\usepackage{caption}
\usepackage{subcaption}
\usepackage{algpseudocode}
\usepackage{comment}
\usepackage{enumerate}
\usepackage[round]{natbib}

\newcommand{\N}{\mathbb{N}}

\newcommand{\R}{\mathbb{R}}

\newcommand{\I}{\mathcal{I}}
\newcommand{\Ip}{\mathcal{I'}}
\newcommand{\Ap}{\mathcal{A'}}
\newcommand{\Ep}{\mathcal{E'}}
\newcommand{\Fp}{\mathcal{F'}}
\newcommand{\comp}{\diamond}
\newcommand{\subgr}{\preceq_{\mathrm{sg}}}

\newtheorem{problem}{Problem}
\newtheorem{lemma}{Lemma}
\newtheorem{proposition}{Proposition}
\newtheorem{theorem}{Theorem}
\newtheorem{definition}{Definition}
\newtheorem{corollary}{Corollary}
\newtheorem{conjecture}{Conjecture}

\usepackage{tikz}

\definecolor{coldegree1}{rgb}{0.62,0.75,0.49}
\definecolor{coldegree2}{rgb}{1.0,1.0,0.0}
\definecolor{coldegree3}{rgb}{1.0,0.0,0.0}
\definecolor{coldegree4}{rgb}{0.0,0.0,1.0}
\definecolor{coldegree5}{rgb}{0.0,0.0,0.0}

\title{The Leaf Function of Penrose P2 Graphs}
\author{Carole Porrier\affiliationmark{1} \and Alain Goupil\affiliationmark{2} \and Alexandre Blondin Massé\affiliationmark{3}}
\affiliation{
  Université Sorbonne Paris Nord, Villetaneuse, France\\
  Université du Québec à Trois-Rivières, Trois-Rivières, Canada\\
  Université du Québec à Montréal, Montréal, Canada}

\begin{document}

\publicationdata{vol. 27:3}{2025}{9}{10.46298/dmtcs.13662}{2024-05-25; 2024-05-25; 2025-06-13}{2025-08-18}
\maketitle

\begin{abstract}
\medskip
We study a graph-theoretic problem in the Penrose P2-graphs which are the dual graphs of Penrose tilings by kites and darts.
Using substitutions, local isomorphism and other properties of Penrose tilings, we construct a family of arbitrarily large induced subtrees of Penrose graphs with the largest possible number of leaves for a given number $n$ of vertices. These subtrees are called fully leafed induced subtrees. 
We denote their number of leaves $L_{P2}(n)$ for any non-negative integer $n$, and the sequence $\left(L_{P2}(n)\right)_{n\in\mathbb{N}}$ is called the leaf function of Penrose P2-graphs.
We  present exact and recursive formulae for $L_{P2}(n)$, as well as an infinite sequence of fully leafed induced subtrees, which are caterpillar graphs.
In particular, our proof relies on the construction of a finite graded poset of $3$-internal-regular subtrees.
\end{abstract}


\section{Introduction}

In graph theory, several problems concern the identification of subgraphs satisfying  particular criteria, such as finding optimal spanning trees \citep{SAT1}, or the {\it snake in the box problem} \citep{snake-in-the-box}.
For instance, in communication networks, it might be cheaper to design a tree topology having as many leaves as possible, since the software and hardware associated to a degree-1 terminal is cheaper than the software and hardware used in the internal terminals that must handle communication and routing \citep{fernandes1998minimal}.
More recently, network topologies with many leaves have  been considered for the coordination of computing clusters \citep{chinnasamy2019minimum}.

\cite{blondin-induced} introduced the following two related problems:
\begin{problem}[\textsc{LIS}]
  Given a simple graph $G$ and two positive integers $i$ and $\ell$, is there an induced subtree of $G$ containing $i$ vertices and $\ell$ leaves?
\end{problem}
\begin{problem}[\textsc{FLIS}]
  Given a simple graph $G$ and a positive integer $i$, what is the largest integer $\ell$ such that there is an induced subtree of $G$ with $i$ vertices and $\ell$ leaves?
\end{problem}

\begin{figure}[h]
    \begin{subfigure}[t]{0.25\textwidth}
      \centering
    \begin{tikzpicture}[scale=0.6, rotate=90, z={(-4mm,-6mm)}, 
        vertex/.style={minimum size=1.8mm, inner sep=0pt, fill=black, circle}, 
        rest/.style={vertex, draw=white, fill=black!20}, 
        subtree/.style={vertex, draw=blue, fill=black, very thick, draw=black}, 
        every edge/.style={draw}, 
        subtreeedge/.style={very thick}]
      \node[rest]    (a) at (-2.3, 1.9) {};
      \node[rest]    (b) at (-2.0, 3.1) {};
      \node[subtree] (c) at (-1.3, 2.4) {};
      \node[subtree] (d) at ( 1.0, 3.1) {};
      \node[subtree] (e) at (-0.3, 3.3) {};
      \node[subtree] (f) at ( 0.2, 1.5) {};
      \node[subtree] (g) at ( 1.3, 1.6) {};
      \node[rest]    (h) at ( 2.4, 0.8) {};
      \node[rest]    (i) at ( 2.4,-0.4) {};
      \node[rest]    (j) at ( 1.8,-2.1) {};
      \node[subtree] (k) at ( 0.8,-1.1) {};
      \node[subtree] (l) at ( 0.0, 0.0) {};
      \node[subtree] (m) at ( 1.0, 0.1) {};
      \node[rest]    (n) at ( 0.6,-2.5) {};
      \node[subtree] (o) at (-0.4,-1.5) {};
      \node[rest]    (p) at (-1.8,-1.4) {};
      \node[subtree] (q) at (-1.0, 1.0) {};
      \node[subtree] (r) at (-1.0,-1.0) {};
      \node[rest]    (s) at (-1.3,-0.1) {};
      \node[rest]    (t) at (-2.5, 0.1) {};
      \path[subtreeedge] (f) edge (c) edge (d) edge (e) edge (g) edge (l);
      \path[subtreeedge] (l) edge (q) edge (r) edge (o) edge (m) edge (k);
      \path (a) edge (b) edge (c) edge (q) edge (t);
      \path (p) edge (t) edge (r) edge (o);
      \path (s) edge (t) edge (q) edge (l) edge (r);
      \path (n) edge (o) edge (k) edge (j);
      \path (j) edge (k) edge (i);
      \path (i) edge (m) edge (h) edge (k);
      \path (h) edge (g) edge (m);
      \path (b) edge (c) edge (e);
    \end{tikzpicture}
    \caption{Finite graph}
    \label{fig:graphe-fini}
    \end{subfigure}
    \hfill
    \begin{subfigure}{0.25\textwidth}
    \vspace{-1cm}
    \centering
    \begin{tikzpicture}[scale=0.6,rotate=90, z={(-4mm,-6mm)}, vertex/.style={minimum size=1.8mm, inner sep=0pt, fill=black, circle}, rest/.style={vertex, draw=white, fill=black!20}, subtree/.style={vertex, draw=blue, fill=black, very thick, draw=black}, degree4/.style={vertex, fill=coldegree4}, degree3/.style={vertex, fill=coldegree3}, degree5/.style={vertex, coldegree5}, every edge/.style={draw}, root/.style={vertex, draw=red, very thick}, fleche/.style={cyan, very thick, ->, shorten >=1mm}, pointe/.style={->}, subtreeedge/.style={very thick}]
      \def\triangleup#1#2#3{\begin{scope}[xshift=#1 cm, yshift=#2 cm]
        \draw[fill=#3] (90:1cm) -- (210:1cm) -- (330:1cm) -- cycle;
        \node[vertex] at (0,0) {};
        \draw (0,0) -- (30:1cm);
        \draw (0,0) -- (150:1cm);
        \draw (0,0) -- (270:1cm);
      \end{scope}}
      \def\triangledown#1#2#3{\begin{scope}[xshift=#1 cm, yshift=#2 cm]
        \draw[fill=#3] (-90:1cm) -- (-210:1cm) -- (-330:1cm) -- cycle;
        \node[vertex] at (0,0) {};
        \draw (0,0) -- (-30:1cm);
        \draw (0,0) -- (-150:1cm);
        \draw (0,0) -- (-270:1cm);
      \end{scope}}
      \def\leafup#1#2#3{\begin{scope}[xshift=#1 cm, yshift=#2 cm]
        \draw[fill=#3] (90:1cm) -- (210:1cm) -- (330:1cm) -- cycle;
        \node[vertex] at (0,0) {};
      \end{scope}}
      \def\leafdown#1#2#3{\begin{scope}[xshift=#1 cm, yshift=#2 cm]
        \draw[fill=#3] (-90:1cm) -- (-210:1cm) -- (-330:1cm) -- cycle;
        \node[vertex] at (0,0) {};
      \end{scope}}
      \begin{scope}[xshift=12.0cm, yshift=-6.8cm, scale=0.7]
        \def\x{0.866}
        \begin{scope}[rotate=90]
        \leafup{-3*\x}{1.5}{coldegree1}
        \leafup{-\x}{1.5}{coldegree1}
        \leafup{\x}{1.5}{coldegree1}
        \leafup{3*\x}{1.5}{coldegree1}
        \leafdown{-2*\x}{-1.0}{coldegree1}
        \leafdown{0}{-1.0}{coldegree1}
        \leafdown{2*\x}{-1.0}{coldegree1}
        \leafup{-4*\x}{0}{coldegree1}
        \leafup{4*\x}{0}{coldegree1}
        \triangledown{-3*\x}{0.5}{coldegree3}
        \triangleup{-2*\x}{0}{coldegree3}
        \triangledown{-\x}{0.5}{coldegree3}
        \triangleup{0}{0}{coldegree3}
        \triangledown{\x}{0.5}{coldegree3}
        \triangleup{2*\x}{0}{coldegree3}
        \triangledown{3*\x}{0.5}{coldegree3}
        \end{scope}
      \end{scope}
    \end{tikzpicture}
    \caption{Polyamond}
    \label{fig:polyamant}
    \end{subfigure}
    \hfill
    \begin{subfigure}[t]{0.28\textwidth}
    \begin{tikzpicture}[scale=0.5]
      \begin{scope}
        \node[scale=1.3] at (0,0) {\includegraphics[scale=0.6]{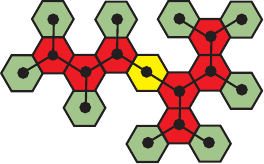}};
      \end{scope}
    \end{tikzpicture}
    \caption{Polyhex}
    \label{fig:Polyhex}
    \end{subfigure}
    
    \begin{subfigure}{0.25\textwidth}
    \centering
    \begin{tikzpicture}[scale=0.6, z={(-4mm,-6mm)}, vertex/.style={minimum size=1.8mm, inner sep=0pt, fill=black, circle}, rest/.style={vertex, draw=white, fill=black!20}, subtree/.style={vertex, draw=blue, fill=black, very thick, draw=black}, degree4/.style={vertex, fill=coldegree4}, degree3/.style={vertex, fill=coldegree3}, degree5/.style={vertex, coldegree5}, every edge/.style={draw}, root/.style={vertex, draw=red, very thick}, fleche/.style={cyan, very thick, ->, shorten >=1mm}, pointe/.style={->}, subtreeedge/.style={very thick}]
      \begin{scope}[xshift=-1cm, yshift=-8cm, scale=0.7]
        \foreach \x/\y in {1/3, 0/1, 0/-1, 4/-1, 1/-2, 3/-2} {
          \draw[ultra thick, fill=coldegree2, rounded corners] (\x,\y) ++ (-0.5,-0.5) rectangle ++ (1,1);
          \node[vertex] at (\x, \y) {};
        }
        \foreach \x/\y in {1/5, 0/4, 2/4, -1/3, -3/2, 3/2, -2/1, 2/1, 4/1, -1/0, 1/0, 3/0, 5/0, 2/-1, -1/-2, 5/-2, 0/-3, 2/-3, 4/-3} {
          \draw[very thick, fill=coldegree1, rounded corners] (\x,\y) ++ (-0.5,-0.5) rectangle ++ (1,1);
          \node[vertex] at (\x, \y) {};
        }
        \foreach \x/\y in {-2/2, 0/2, 2/2} {
          \draw[very thick, fill=coldegree3, rounded corners] (\x,\y) ++ (-0.5,-0.5) rectangle ++ (1,1);
          \node[vertex] at (\x, \y) {};
          \draw[thick] (\x,\y) -- ++ (1,0);
          \draw[thick] (\x,\y) -- ++ (-1,0);
          \draw[thick] (\x,\y) -- ++ (0,-1);
        }
        \foreach \x/\y in {-1/2, 1/2} {
          \draw[very thick, fill=coldegree3, rounded corners] (\x,\y) ++ (-0.5,-0.5) rectangle ++ (1,1);
          \node[vertex] at (\x, \y) {};
          \draw[thick] (\x,\y) -- ++ (1,0);
          \draw[thick] (\x,\y) -- ++ (-1,0);
          \draw[thick] (\x,\y) -- ++ (0,1);
        }
        \foreach \x/\y in {0/0, 4/0, 0/-2, 2/-2, 4/-2, 1/4} {
          \draw[very thick, fill=coldegree4, rounded corners] (\x,\y) ++ (-0.5,-0.5) rectangle ++ (1,1);
          \node[vertex] at (\x, \y) {};
          \draw[thick] (\x,\y) -- ++ (1,0);
          \draw[thick] (\x,\y) -- ++ (-1,0);
          \draw[thick] (\x,\y) -- ++ (0,1);
          \draw[thick] (\x,\y) -- ++ (0,-1);
        }
      \end{scope}
    \end{tikzpicture}
    \caption{Polyomino}
    \label{fig:Polyomino}
    \end{subfigure}
    \hfill
    \begin{subfigure}{0.65\textwidth}
    \centering
    \begin{tikzpicture}[scale=0.85, z={(-4mm,-6mm)}, vertex/.style={minimum size=1.8mm, inner sep=0pt, fill=black, circle}, rest/.style={vertex, draw=white, fill=black!20}, subtree/.style={vertex, draw=blue, fill=black, very thick, draw=black}, degree4/.style={vertex, fill=coldegree4}, degree3/.style={vertex, fill=coldegree3}, degree5/.style={vertex, coldegree5}, every edge/.style={draw}, root/.style={vertex, draw=red, very thick}, fleche/.style={cyan, very thick, ->, shorten >=1mm}, pointe/.style={->}, subtreeedge/.style={very thick}]
  \def\a{0.5}
  \def\zscale{1.4}
  \begin{scope}[xshift=8cm, xscale=0.6, yscale=0.6]
      \draw (-6,0,0) -- (8,0,0);
      \draw (-3,-3,0) -- (-3,3,0);
      \draw (0,0,-3) -- (0,0,3);
      \draw (-3,3,-1) -- (-3,3,1);
      \draw (-3,-3,-1) -- (-3,-3,1);
      \draw (0,-1,3) -- (0,1,3);
      \draw (0,-1,-3) -- (0,1,-3);
      \draw (-6,0,-1) -- (-6,0,1);
      \draw (8,0,-1) -- (8,0,1);
      \foreach \x/\y/\z in {-5/0/0, -3/3/1, -3/3/-1, -3/-3/1, -3/-3/-1, 0/1/-3, 0/-1/-3, 0/1/3, 0/-1/3, -3/2/0, -3/-2/0, 0/0/-2, 0/0/2, 0/0/0} {
        \draw (\x,\y,\z) ++ (-\a,0,0) -- ++ (2*\a,0,0);
      }
      \foreach \x/\y/\z in {-6/0/1, -6/0/-1, -5/0/0, -1/0/0, 0/0/-1, 0/0/1, 0/-1/-3, 0/1/-3, 0/-1/3, 0/1/3, 3/0/0, 5/0/0, 7/0/0, 8/0/1, 8/0/-1} {
        \draw (\x,\y,\z) ++ (0,-\a,0) -- ++ (0,2*\a,0);
      }
      \foreach \x/\y/\z in {-6/0/1, -6/0/-1, -4/0/0, -2/0/0, -3/1/0, -3/-1/0, -3/-3/-1, -3/-3/1, -3/3/-1, -3/3/1, 2/0/0, 4/0/0, 6/0/0, 8/0/1, 8/0/-1} {
        \draw (\x,\y,\z) ++ (0,0,-\zscale*\a) -- ++ (0,0,2*\zscale*\a);
      }
      \foreach \x/\y/\z in {-6/0/1, -6/0/-1} {
        \draw (\x,\y,\z) -- ++ (-\a,0,0);
      }
      \foreach \x/\y/\z in {8/0/1, 8/0/-1} {
        \draw (\x,\y,\z) -- ++ (\a,0,0);
      }
      \foreach \x/\y/\z in {-3/-3/1, -3/-3/-1, 1/0/0} {
        \draw (\x,\y,\z) -- ++ (0,-\a,0);
      }
      \foreach \x/\y/\z in {-3/3/1, -3/3/-1} {
        \draw (\x,\y,\z) -- ++ (0,\a,0);
      }
      \foreach \x/\y/\z in {0/-1/-3, 0/1/-3} {
        \draw (\x,\y,\z) -- ++ (0,0,-\zscale*\a);
      }
      \foreach \x/\y/\z in {0/-1/3, 0/1/3} {
        \draw (\x,\y,\z) -- ++ (0,0,\zscale*\a);
      }
    \foreach \x/\y/\z in {0/0/0, -1/0/0, -2/0/0, -3/0/0, -4/0/0, -5/0/0, -3/1/0, -3/2/0, -3/-1/0, -3/-2/0, 0/0/1, 0/0/2, 0/0/-1, 0/0/-2, 2/0/0, 3/0/0, 4/0/0, 5/0/0, 6/0/0, 7/0/0} {
      \node[degree4] (\x \y \z)  at (\x, \y, \z)   {};
    }
    \foreach \x/\y/\z in {-6/0/0, -3/3/0, -3/-3/0, 0/0/3, 0/0/-3, 1/0/0, 8/0/0} {
      \node[degree3] (\x \y \z)  at (\x, \y, \z)   {};
    }
    \foreach \x/\y/\z in {-6/0/1, -6/0/-1, -3/3/1, -3/3/-1, -3/-3/1, -3/-3/-1, 0/1/3, 0/-1/3, 0/1/-3, 0/-1/-3, 8/0/1, 8/0/-1} {
      \node[degree5] (\x \y \z)  at (\x, \y, \z)   {};
    }
  \end{scope}
\end{tikzpicture}
    \caption{Structure of a polycube}
    \label{fig:polycube}
    \end{subfigure}
    \caption{Fully leafed induced subtrees, (a) in a finite graph, (b) in the triangular lattice, (c) in the hexagonal lattice, (d) in the square lattice, and (e) in the cubic lattice. 
    Vertices are colored according to their degree in the induced subtree (except for the finite graph).
    Subtrees (a) and (b) are caterpillars, subtrees (c), (d) and (e) are not.}
    \label{fig:FLIS1}
\end{figure}
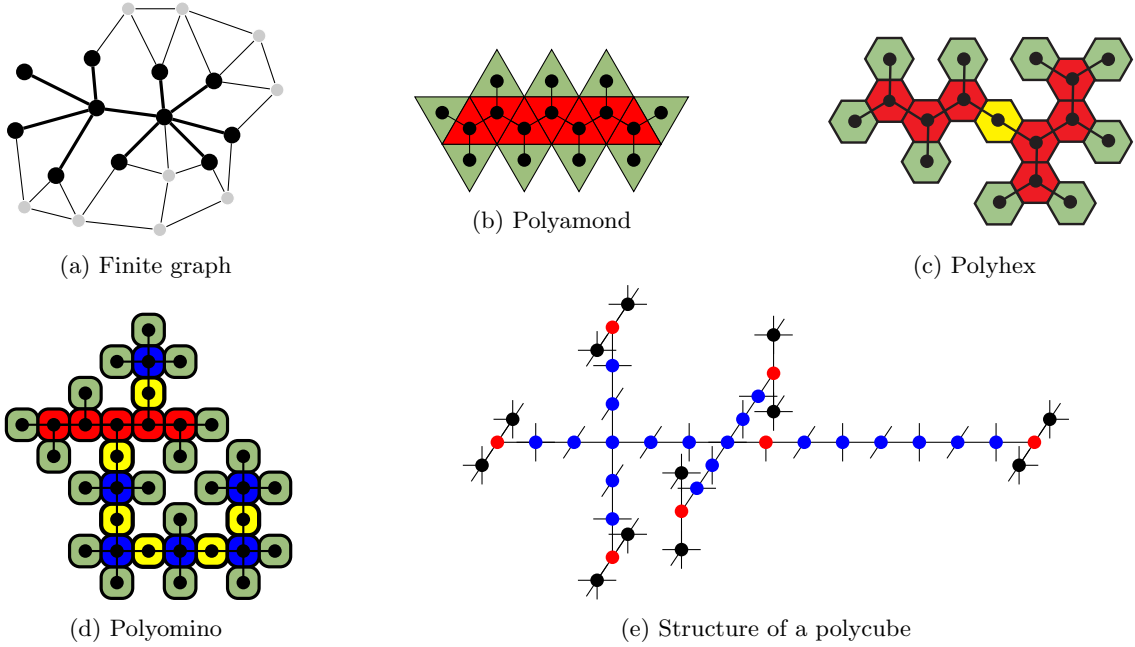

Induced subtrees maximizing their number of leaves with respect to their number of vertices are called {\it fully leafed}.
Some examples are illustrated in Figure \ref{fig:FLIS1}.
As it is often the case in graph theoretical problems, it was shown that the LIS problem is NP-complete, so that the FLIS problem is NP-hard \citep{blondin-induced}.
The same authors proposed a polynomial time algorithm using dynamic programming for the special case where $G$ is a tree.
They also provided a general exponential branch-and-bound algorithm that solves the problem for any graph.

In the same article, another line of investigation was explored on infinite graphs.
Indeed, given some infinite simple graph $G$, it was proposed to study the function $L_G(i)$, called the \emph{leaf function} of $G$, which associates with every  integer $i\geq 0$, the maximal number of leaves that an induced subtree of $G$ of order $i$  can have.
The leaf functions of the square grid, the triangular grid, the hexagonal grid and the cubic grid have thus been described by \cite{blondin-poly} and
 some special classes of fully leafed induced subtrees have been enumerated \citep{blondin-saturated}.
Figure \ref{fig:FLIS1} illustrates some fully leafed induced subtrees in these infinite graphs, which are actually dual graphs of tilings.

Here we focus on Penrose tilings, which are aperiodic.
The discovery of quasicrystals in the 1980s provoked a small revolution in classical crystallography, and Penrose tilings serve as models for them \citep{Senechal1995}. 
Since their discovery, extensive studies were conducted on the properties of quasicrystals \citep{Vedmedenko01052008,CWS2017}. 
Fully leafed induced subtrees of any 3D network satisfy the property of having maximal area for a given volume and in 2D networks, they maximize perimeter with respect to area. This physical property is important in the interaction between the surface of a solid and its environment. Quasicrystal surfaces are emerging as interesting materials for adsorption and reactions on surfaces \citep{McGrath_2010,DUBOIS2000,KWEON2022117657}, low friction machine parts and nonstick coatings \citep{DUBOIS2000,Diehl_2008}. 

Moreover, the Potts model has been used as a simple model of protein folding \citep{LHTW1996,Garel_1988,desilva2018}.  Compared to the typical approximation that the protein folds on a square lattice \citep{ChanDill1989,bodroza2013,desilva2018}, aperiodic networks arguably define a better approximation to the real problem  because aperiodic networks have a larger range of bond angles, providing a better approximation of a continuous rotational symmetry than any 2D regular lattice.

\begin{figure}
    \centering
    \input{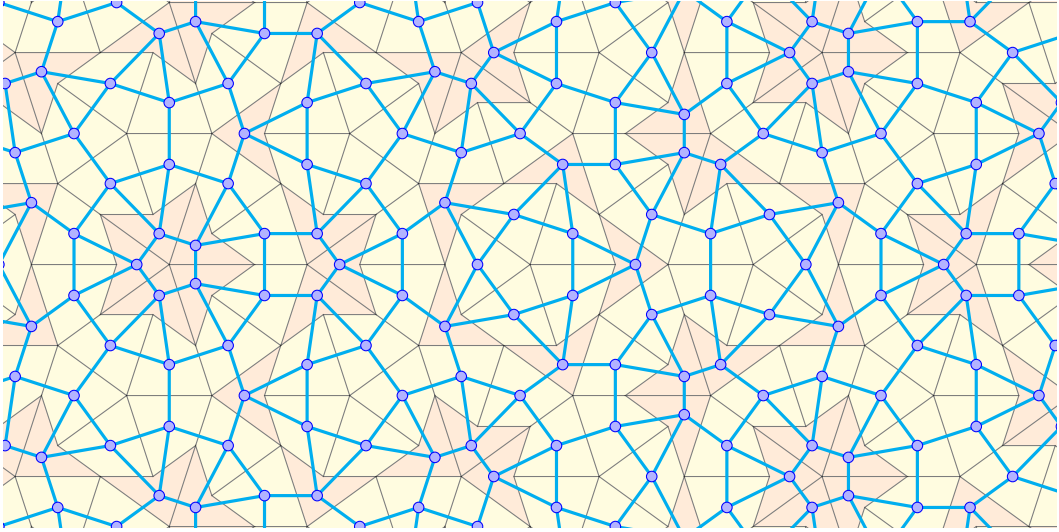}
    \caption{An excerpt of a P2 tiling, also called \emph{kites and darts tiling}, with its adjacency graph superimposed. The kites and darts are respectively colored yellow and orange. The graph vertices are depicted with blue dots and the edges in cyan.}

    \label{fig:penrose-graph}
\end{figure}

In this paper, we investigate the leaf function $L_{P{2}}$ of Penrose tilings of type $P{2}$, also known as \emph{kites and darts} tilings. 
The adjacency graph of a tiling of type $P2$ is called a {\it P2-graph} or a {\it Penrose tree} (resp. {\it Penrose caterpillar}) when it is a tree (resp. a caterpillar).
A P2 tiling with its graph superimposed is illustrated in Figure \ref{fig:penrose-graph}.
We prove that fully leafed induced subtrees in P2 tilings are essentially caterpillars and we present an explicit procedure for their construction in P2. We believe that these caterpillars and their corresponding paths in P2 are of interest in the study of quasicrystals and possibly in protein folding. 
Our main contributions are the construction of arbitrarily large fully leafed induced subtrees and also the computation of the leaf function for Penrose tilings by kites and darts, summarized in the following theorem:

\begin{theorem}\label{thm1}
The leaf function $L_{P2}$ of P2-graphs, defined for all $n\in\N$, is
$$
L_{P2}(n) = 
\begin{cases}
    0 & \text{if } ~ n\in \{ 0, 1\}\\
    \left\lfloor n/2 \right\rfloor +1 & \text{if } ~2\leq n\leq 18,\\
    L_{P2}(n-17)+8 & \text{if } ~n\geq 19.
\end{cases}
$$
\end{theorem}

The remaining sections of this article are divided as follows.
Section 2 introduces the definitions and notations on tilings and graphs.
In Sections 3 and 4, we prove that the function given in Theorem \ref{thm1} is respectively an upper bound and a lower bound for the leaf function of P2 graphs.
Finally, in Section 5 we give a closed formula for $L_{P2}$, we discuss perspectives and we conclude with two conjectures.


\section{Leaf functions and Penrose trees}

We start by recalling usual graph theoretical definitions and notation.

Let $G = (V, E)$ be a simple undirected graph, finite or infinite, with set of vertices $V$ and set of edges $E$. 
Let $u\in V$ and $U \subseteq V$.
Then the set of {\it neighbors} of $u$ in $G$ is denoted $N_G(u)$, which is naturally extended to $U$ by defining its {\it neighborhood} $N_G(U) = \{N_G(u') \mid u' \in U\}$.
The {\it subgraph of $G$ induced by $U$} is $G[U] = (U, E \cap \mathcal{P}_2(U))$, where $\mathcal{P}_2(U)$ is the set of 2-element subsets of $U$. 
The {\it degree} of a vertex is its number of neighbors.

A {\it tree} is a connected and acyclic graph.
The set of all subgraphs of $G$ can be equipped with a partial order $\subgr$: 
if $H_1,H_2$ are two subgraphs of a graph $G$, we write $H_1\subgr\;H_2$ when $H_1$ is a subgraph of $H_2$ up to isomorphism, i.e., if there exist two graphs $\Tilde{H}_1,\Tilde{H}_2$ such that $H_1$ is isomorphic to $\Tilde{H}_1$,  $H_2$ is isomorphic to $\Tilde{H}_2$ and $\Tilde{H}_1$ is a subgraph of $\Tilde{H}_2$.

We focus on \emph{induced subtrees}, i.e., the induced subgraphs of $G$ that are trees, and in their {\it leaves}, i.e., their vertices of degree $1$.
Given an induced subtree $I$, the induced subtree obtained from $I$ by deleting all its leaves is called the \emph{derived subtree} of $I$ and is denoted by $I'$.
An induced subtree $I$ is a \emph{chain} if it is a {\it path graph}, i.e., a tree whose vertices have degree at most $2$.
An induced subtree $I$ is a \emph{caterpillar} if $I'$ is a path graph.
We denote by $n(I)$, the total number of vertices of an induced subtree $I$ and by $n_{d}(I)$, the number of vertices of degree $d$ in $I$.

Let $\mathcal{I}_G(n)$ be the set of all induced subtrees of $G$ of order $n$, up to isomorphism.
Then the {\it leaf function} of $G$ is defined  by 
$$L_G(n)= \max\{n_1(I) : I\in\mathcal{I}_G(n)\},$$
for all $n\in\N$. 
An induced subtree $I$ of $G$ of order $n$ is said to be {\it fully leafed} if $n_1(I) = L_G(n)$.
For any positive integer $d$, an induced subtree $I$ is said to be {\it $d$-internal-regular}  when all the vertices of $I'$ are of degree $d$ in $I$. 
In the following section, we focus on $3$-internal-regular subtrees of kites and darts Penrose tilings.


The graphs that we are interested in here are defined on tilings.
A \emph{tiling} of the plane is a covering of $\R^2$ by countably many closed sets called \emph{tiles} with no gaps and no overlaps, except for their boundary points which form the vertices and edges of the tiling.
Two tiles are said to be \emph{adjacent} if they share an edge.
The vertices and edges of a tiling thus form a planar (undirected) graph embedded in the plane, whose dual is the {\it adjacency graph} $G$ of the tiling:
the vertices of $G$ are the tiles of the tiling and the set $E$ of edges is given by the adjacency relation between the tiles.
In that sense, polyominoes, polyamonds and polyhexes like those in Figure \ref{fig:FLIS1} are \emph{patches}, i.e., connected sets of tiles, but can also be seen as induced subtrees of (the duals of) the three regular tilings of the plane -- respectively the square, triangular and hexagonal lattices.

\begin{figure}[h]
    \begin{subfigure}[b]{0.18\textwidth}
        \includegraphics[scale=0.35]{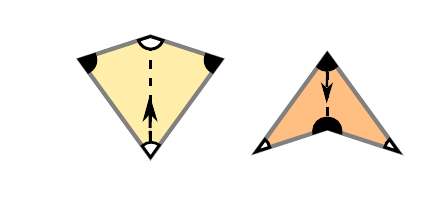}
        \caption*{\emph{kite} \hspace{0.2cm} \emph{dart}}
        \label{kd}
    \end{subfigure}
    \hfill
    \begin{subfigure}[b]{0.11\textwidth}
        \centering
        \includegraphics[scale=0.35]{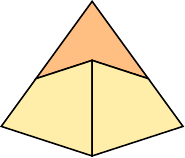}
        \caption*{\emph{Ace}}
        \label{ace}
    \end{subfigure}
    \hfill
    \begin{subfigure}[b]{0.11\textwidth}
        \centering
        \includegraphics[scale=0.35]{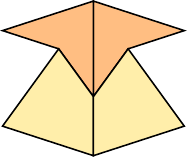}
        \caption*{\emph{Deuce}}
        \label{deuce}
    \end{subfigure}
    \hfill
    \begin{subfigure}[b]{0.11\textwidth}
        \centering
        \includegraphics[scale=0.35]{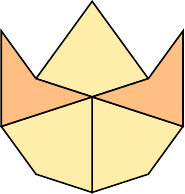}
        \caption*{\emph{Jack}}
        \label{jack}
    \end{subfigure}
    \hfill
    \begin{subfigure}[b]{0.11\textwidth}
        \centering
        \includegraphics[scale=0.35]{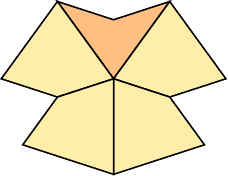}
        \caption*{\emph{Queen}}
        \label{queen}
    \end{subfigure}
    \hfill
    \begin{subfigure}[b]{0.11\textwidth}
        \centering
        \includegraphics[scale=0.35]{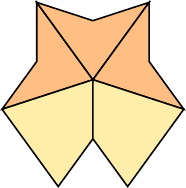}
        \caption*{\emph{King}}
        \label{king}
    \end{subfigure}
    \hfill
    \begin{subfigure}[b]{0.11\textwidth}
        \centering
        \includegraphics[scale=0.35]{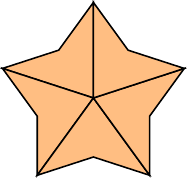}
        \caption*{\emph{Star}}
        \label{star}
    \end{subfigure}
    \hfill
    \begin{subfigure}[b]{0.11\textwidth}
        \centering
        \includegraphics[scale=0.35]{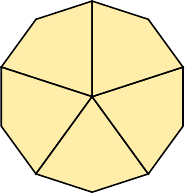}
        \caption*{\emph{Sun}}
        \label{sun}
    \end{subfigure}
    \caption{The kite and dart, followed by the seven vertex configurations in a Penrose P2 tiling. The colored vertices on the kite and dart on the left give the assembly rules: the colors have to match as one assembles the tiles, so that there is either a black or a white circle at each vertex. No other vertex configuration than those above is allowed -- in particular, no rhombus can be formed. All patches are considered up to isometry.}
    \label{vertex-config}
\end{figure}

The main objects of this article are Penrose tilings by kites and darts, classified as P2 in Grünbaum and Shephard's reference book  (\citeyear{grunbaum2016}).
For convenience, we call the dual graph of a Penrose tiling a \emph{P2-graph}, and an induced subgraph of a P2-graph is called a \emph{Penrose tree} (resp. \emph{Penrose caterpillar}) when it is a tree (resp. a caterpillar).
A $P2$ tiling is composed of two kinds of tiles called \emph{kites} and \emph{darts}, as shown on the left of Figure \ref{vertex-config}.
With the P2 assembly rules, there are seven possible \emph{vertex configurations}, i.e., seven ways of arranging the tiles around a vertex of the plane (also in Figure \ref{vertex-config}). 

\begin{figure*}[h]
    \centering
    \begin{subfigure}[b]{0.08\textwidth}
        \centering
        \includegraphics[scale=0.3]{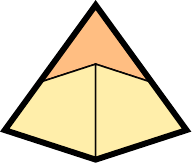}
        \caption*{\centering Ace}
    \end{subfigure}
    \hfill
    \begin{subfigure}[b]{0.1\textwidth}
        \centering
        \includegraphics[scale=0.3]{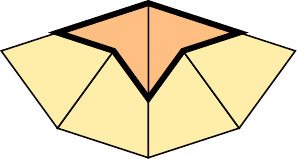}
        \caption*{\centering Deuce}
    \end{subfigure}
    \hfill
    \begin{subfigure}[b]{0.12\textwidth}
        \centering
        \includegraphics[scale=0.3]{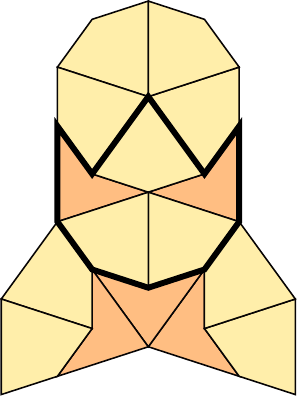}
        \caption*{\centering Jack}
    \end{subfigure}
    \hfill
    \begin{subfigure}[b]{0.18\textwidth}
        \centering
        \includegraphics[scale=0.3]{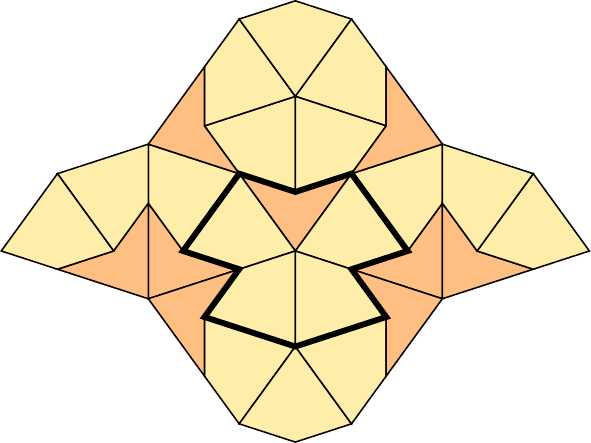}
        \caption*{\centering Queen}
    \end{subfigure}
    \hfill
    \begin{subfigure}[b]{0.2\textwidth}
        \centering
        \includegraphics[scale=0.3]{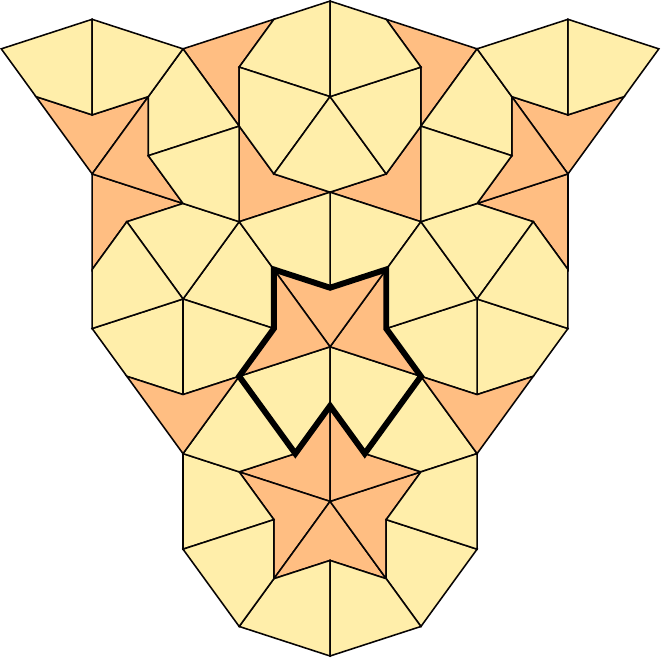}
        \caption*{\centering King}
    \end{subfigure}
    \hfill
    \begin{subfigure}[b]{0.1\textwidth}
        \centering
        \includegraphics[scale=0.3]{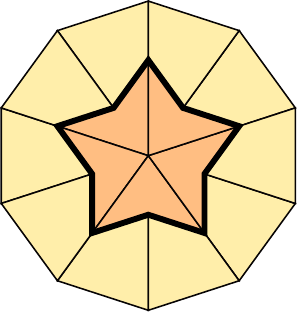}
        \caption*{\centering Star}
    \end{subfigure}
    \hfill
    \begin{subfigure}[b]{0.08\textwidth}
        \centering
        \includegraphics[scale=0.3]{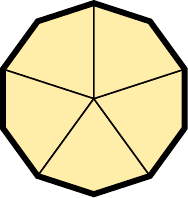}
        \caption*{\centering Sun}
    \end{subfigure}
    \caption{The seven kingdoms of kites and darts tilings.}
    \label{fig:kingdoms-P2}
\end{figure*}

Each vertex configuration \emph{forces} a whole set of tiles around it which John H. Conway called an \emph{empire}  \citep{gardner1977,grunbaum2016}.
This means that whenever a given vertex configuration appears in a P2 tiling, then all the tiles forming the empire are placed in a unique possible way with respect to that vertex.
Some empires are not connected and even contain infinitely many tiles, yet it is sufficient here to consider, for each vertex configuration, the largest connected component included in its empire. We call this component its \emph{kingdom} (or \emph{local empire}).
The seven kingdoms of kites and darts tilings are shown in Figure \ref{fig:kingdoms-P2}.

\begin{figure}[ht!]
    \centering
    \includegraphics[width=0.85\textwidth]{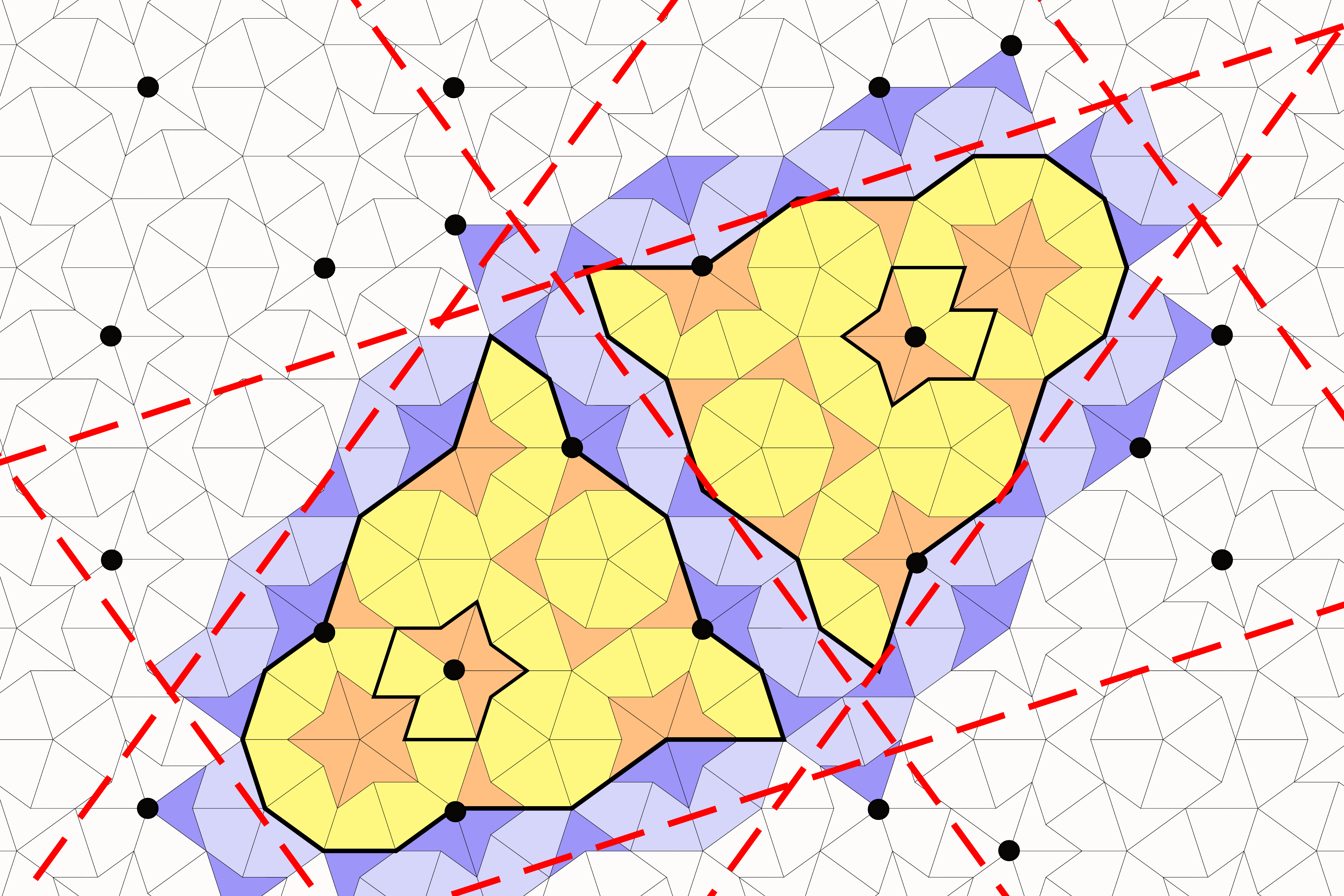}
    \caption{Each black point corresponds to an occurrence of a King, with the same kingdom around. Two occurrences of the King's kingdom are colored in yellow and orange, in different orientations and separated by a (finite) Conway worm colored in blue, as well as the other finite worms surrounding the highlighted kingdoms. Each worm could be flipped around a central axis parallel to the \emph{Ammann bar} (red line) crossing the worm -- the bar being flipped around the same axis in that case, for instance touching the kingdom on the left instead of that on the right.}
    \label{fig:kingdom-worms}
\end{figure}
As shown in Figure \ref{fig:kingdom-worms}, kingdoms are surrounded by \emph{Conway worms}\footnote{We sometimes just say ``worm'', omitting ``Conway''.}, consisting of long and short \emph{bowties} placed side by side. Those are described in detail by \cite{grunbaum2016}, as well as Ammann bars whose main properties are recalled in Section \ref{optimal-paths}. 

A vertex configuration, a kingdom or a Conway worm, is a special case of a patch.
All Penrose tilings of a given type are \emph{locally isomorphic}, meaning that any patch which appears in one also appears in all the others. 
In other words, all kites and darts Penrose tilings are in the same local isomorphism class (\emph{LI class}).
Local isomorphism and the following notions are explained in more detail in several books and surveys \citep{GBS1994,baake1999guide,Robinson2004,baake_grimm_2013}.

A tiling $T_2$ is \emph{locally derivable} from a tiling $T_1$ if there exists a fixed radius $r$ such that $T_{2}$ can be obtained from $T_{1}$ by a \emph{local mapping}, that is a transformation which maps any patch of radius $r$ in $T_1$ to a given patch in $T_2$, which is always the same up to translation.
If $T_1$ is also locally derivable from $T_2$, then $T_1$ and $T_2$ are \emph{mutually locally derivable} (\emph{MLD}).
The MLD property holds for Penrose tilings of different types (i.e., of different LI classes): any Penrose tiling of one type is MLD with a Penrose tiling of each type, including some of the same LI class, and there is a bijection between the sets of tilings of these two types.
The Penrose tilings are therefore all in the same MLD class, along with a few other sets of tilings including the so-called \emph{HBS tilings}, introduced in Section 4.

Some cases of local derivability are fairly easy to see.
For example, when a tiling $T_1$ is a \emph{composition} of a tiling $T_2$, i.e., when each tile of $T_1$ is a union of tiles of $T_2$ (the pattern being the same for all copies of the same tile). 
In this case, $T_2$ is a \emph{decomposition} of $T_1$.
These notions can be relaxed in the \emph{imperfect} version, where each tile of $T_1$ is \textsl{covered} by a union of tiles of $T_2$.
For Penrose tilings and their derivatives, imperfect composition and decomposition rules are available to transform a tiling of a given type into another tiling of the same or another LI class.
In particular, HBS tiles are a composition of darts and half-kites (Figure \ref{fig:star-tileset}).

An \emph{$\alpha$-decomposition} is a decomposition in which the new tiles have the same shapes as the original tiles but on a $\alpha:1$ scale.
An \emph{inflation} is a decomposition followed by a scaling, so that the new tiles have the same shapes and sizes as the original tiles, but form a larger pattern.
For the Penrose MLD class, the scaling is of ratio $\varphi:1$, so the inflation is a $\varphi$-decomposition, where $\varphi=(1+\sqrt{5})/2$ is the golden ratio.


\section{An upper bound on the Leaf Function}

In this section, we prove that the expression given in Theorem 1 is indeed an upper bound on the function $L_{P2}$, i.e.,

\begin{proposition}\label{prop:borne-sup}
For any non negative integer $n$,
\begin{equation}\label{ineq:upper-bound}    
L_{P2}(n) \leq 
\begin{cases}
    0 & \text{if } ~0\leq n\leq 1,\\
    \left\lfloor n/2 \right\rfloor +1 & \text{if } ~2\leq n\leq 18,\\
    L_{P2}(n-17)+8 & \text{if } ~n\geq 19.
\end{cases}
\end{equation}
\end{proposition}

In order to prove Proposition \ref{prop:borne-sup}, we proceed in five steps:
\begin{enumerate}
    \item We first notice that in any induced subtree of a P2-graph $G$, each internal tile has degree at most $3$ (Lemma \ref{lem:deg-at-most-3});
    \item We provide a relation between $n_1(I)$ and $n_3(I)$ for any induced subtree $I$ of $G$ (Lemma \ref{lem:n1n3});
    \item  We provide a loose upper bound on $L_{P2}$ (Lemma \ref{lem:loose-upper-bound});
    \item We show that any $3$-internal-regular induced subtree of $G$ has at most $8$ internal tiles (Lemma \ref{lem:3-regular});
    \item Using a minimal counter-example argument, we show that Inequation \eqref{ineq:upper-bound} holds.
\end{enumerate}

The first step is straightforward.

\begin{lemma}[\cite{PB2020}]\label{lem:deg-at-most-3}
Let $I$ be an induced subtree of a P2-graph. Then for any tile $t$ in $I$, $\deg_{I}(t) \leq 3$.
\end{lemma}

\begin{figure}[h]
    \centering
    \includegraphics[scale=0.4]{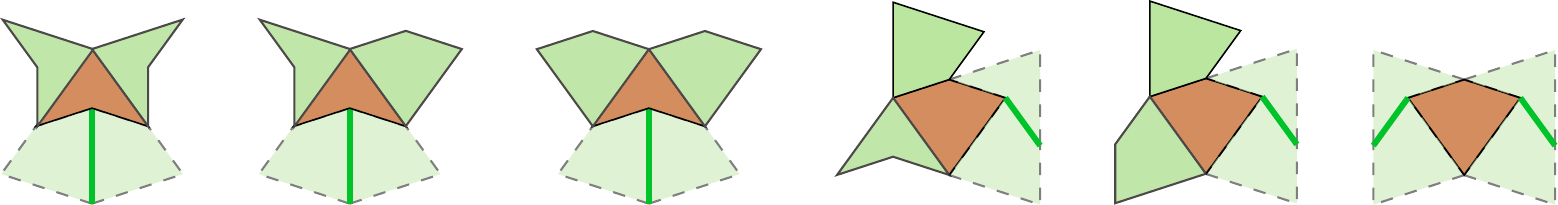}
    \caption{The six possible neighborhoods, up to isometry, of a P2 tile. A green edge between two lighter tiles indicates that they cannot both be in the same induced subtree with the brown tile, since the three tiles  form a cycle.}\label{fig:neigbhorhood}
\end{figure}

\begin{figure}
    \includegraphics[scale=0.23]{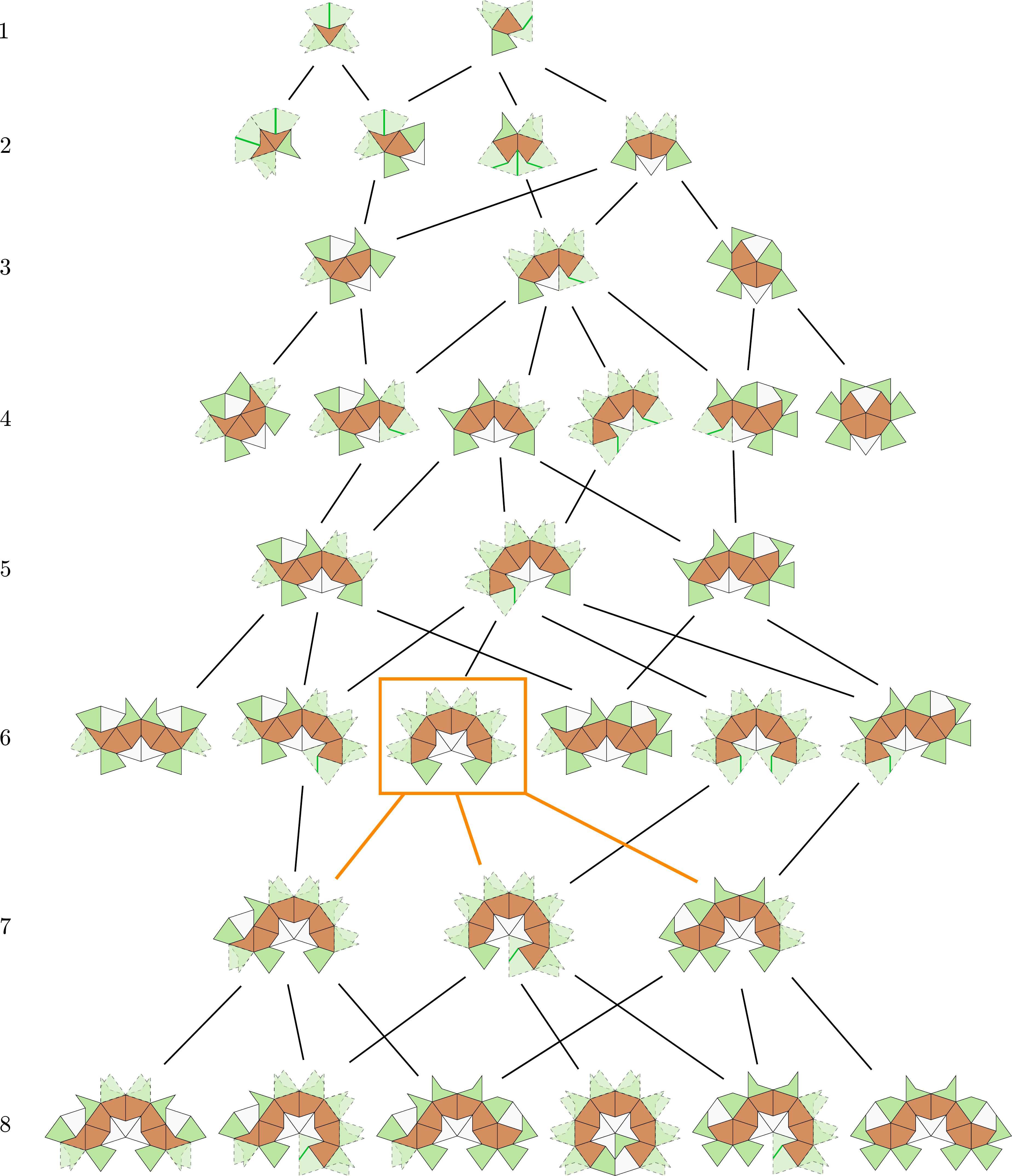}
    \caption{The graded poset $\Ip$ of derived path of 3-internal-regular caterpillars of length 1 to 8, up to isometry. 
             For each caterpillar, the derived path is brown, its adjacent tiles which cannot be in an induced subtree are white and the leaves are green.
             Translucent leaves with dashed contour indicate two possibilities: two tiles are superimposed or adjacent through a green edge, and only one of them is in the induced subtree.
           $\Ip(n)$ is the set of brown paths in row $n$.
             $\Ip(6)$ contains the special caterpillar $C_{14}$ (composed of 14 tiles), in the orange box: its derived path is included in all the caterpillars below.}
    \label{fig:poset}
\end{figure}

\begin{proof}
\cite{PB2020} have proved this result, but for sake of completeness, we include it here.
It suffices to consider all possible neighborhoods of a kite and a dart in a P2 tiling, shown in Figure \ref{fig:neigbhorhood} (which details the first row of Figure \ref{fig:poset}, plus one additional configuration).
In each case, there is at least one pair of adjacent neighbors among the four neighbors of a tile, so that the $I$-degree of $t$ is at most $3$ for any induced subtree $I$.
\end{proof}

The two following lemmas are consequences of elementary graph theory arguments and of Lemma \ref{lem:deg-at-most-3}.

\begin{lemma}\label{lem:n1n3}
    Let $I$ be an induced subtree of a P2-graph, containing $n$ tiles. Then
    \begin{equation}\label{eq:n1n3}
        n_1(I) = n_3(I) + 2.
    \end{equation}
\end{lemma}

\begin{proof}
    For $d \in\{ 1,2,3\}$, let $n_d = n_d(I)$.
    From Lemma \ref{lem:deg-at-most-3}, we have $n = n_1 + n_2 + n_3$ i.e.,
    \begin{equation}\label{eq:n}
      n_2 = n - n_1 - n_3.
    \end{equation}
    Since $I$ is a tree, it has $n - 1$ edges so by the handshaking lemma, $2(n - 1) = n_1 + 2n_2 + 3n_3$ i.e.,
    \begin{equation}\label{eq:handshaking}
      2n_2 = 2(n - 1) - n_1 - 3n_3.
    \end{equation}
    Combining Equations \eqref{eq:n} and \eqref{eq:handshaking} yields the result.
\end{proof}

\begin{lemma}\label{lem:loose-upper-bound}
    For all integers $n, k$ such that $1 \leq k < n$,
    \begin{equation}\label{eq:loose-upper-bound}
    L_{P_{2}}(n) \leq L_{P_{2}}(n - k) + \lceil k / 2 \rceil.
    \end{equation}
\end{lemma}

\begin{proof}
    First, we prove that $L_{P_{2}}(n) \leq L_{P_{2}}(n - 1) + 1$.
    Arguing by contradiction, assume that there exists some integer $n$ such that $L_{P_{2}}(n) > L_{P_{2}}(n - 1) + 1$.
    Let $I$ be a fully leafed induced subtree of size $n$, i.e., such that $n_1(I)=L_{P_{2}}(n)$, and $J$ a subtree obtained by removing any one leaf from $I$.
    By construction of $J$, we have $n_1(J)\geq n_1(I)-1$ and $J$ has size $n-1$, so $n_1(J)\leq L_{P_{2}}(n-1)$.
    Hence $L_{P_{2}}(n - 1)\geq n_1(J)\geq n_1(I)-1 = L_{P_{2}}(n)-1$ in contradiction with the assumption.
    Next, we show 
    \begin{align}\label{eqL(n-2)}
    L_{P_{2}}(n) \leq L_{P_{2}}(n - 2) + 1.
    \end{align}
    From the previous paragraph, $L_{P_{2}}(n-1) \leq L_{P_{2}}(n - 2) + 1$.
    Thus in case $L_{P_{2}}(n) = L_{P_{2}}(n - 1)$, we have $L_{P_{2}}(n) \leq L_{P_{2}}(n - 2) + 1$.
    Otherwise, by Lemma 1 and Proposition 1 of \cite{PB2020}, if $L_{P_{2}}(n) > L_{P_{2}}(n - 1)$ then $L_{P_{2}}(n+1) = L_{P_{2}}(n)$.
    As a consequence, if $L_{P_{2}}(n) > L_{P_{2}}(n - 1)$ then $L_{P_{2}}(n-1) = L_{P_{2}}(n-2)$, which also implies $L_{P_{2}}(n) \leq L_{P_{2}}(n - 2) + 1$.\\
    We are now ready to prove the general case by induction on $k$.

    \textsc{Basis}. The cases $k = 1$ and $k = 2$ have been proved in the previous paragraphs.

    \textsc{Induction hypothesis}. Assume that for some integer $k \geq 3$ and all naturals $k' < k$, one has $L_{P_{2}}(n) \leq L_{P_{2}}(n - k') + \lceil k' / 2 \rceil$.

    \textsc{Induction}. Then
    $$\begin{array}{rcll}
        L_{P_{2}}(n - k) + \lceil k / 2 \rceil
            & = & L_{P_{2}}((n - 2) - (k - 2)) + \lceil k / 2 \rceil \\
            & = & L_{P_{2}}((n - 2) - (k - 2)) +~\lceil (k - 2) / 2 \rceil + 1 \\
            & \geq & L_{P_{2}}(n - 2) + 1 & \mbox{by induction hypothesis} \\
            & \geq & L_{P_{2}}(n), & \mbox{by Equation }\eqref{eqL(n-2)}
    \end{array}$$
    concluding the proof.
\end{proof}

We now show that the set $\I$ of all $3$-internal-regular Penrose trees is finite by exhibiting the complete set $\Ip$ of their derived subtrees, up to isomorphism, in the form of a graded poset $(\Ip,\subgr)$, where $\subgr$ is the partial order relation defined as follows:
if $H_1,H_2$ are two subgraphs of a graph $G$, we write $H_1\subgr\;H_2$ when $H_1$ is a subgraph of $H_2$ up to isomorphism, i.e., if there exist two graphs $\Tilde{H}_1,\Tilde{H}_2$ such that $H_1$ is isomorphic to $\Tilde{H}_1$,  $H_2$ is isomorphic to $\Tilde{H}_2$ and $\Tilde{H}_1$ is a subgraph of $\Tilde{H}_2$.

\begin{lemma}\label{lem:3-regular}
Let $\I$ be the set of all $3$-internal-regular Penrose trees and $\Ip = \{I' \mid I \in \I\}$, i.e., $I' \in \Ip$ if and only if $I'$ is the derived tree of some $3$-internal-regular Penrose tree.
Then $(\Ip,\subgr)$ is the finite graded poset illustrated in Figure \ref{fig:poset}. In particular, $\I$ is finite.
\end{lemma}

\begin{proof}
    Let us first introduce some useful sets. 
    For each positive integer $n$, let $\Ip(n)$ be the set of all induced subtrees $I'$ such that $I'$ has exactly $n$ tiles and there exists a $3$-internal-regular induced subtree whose derived tree is $I'$.
    We say that an induced subtree $E$ is an \emph{extension} of a subtree $I$ if $n(E) = n(I) + 1$ and $I$ is a subtree of $E$.
    For each positive integer $n$, let $\Ep(n)$ be the set of all extensions $E$ of subtrees $I\in \Ip(n-1)$, 
    and let $\Ap(n)$ be the subset of $\Ep(n)$ consisting of induced subtrees $A$ such that each subtree of $A$ of size $n' < n$ belongs to $\Ip(n')$.
    To make the argument more intuitive, we call an element $A$ of $\Ap(n)$ \emph{admissible} in the sense that all its proper subtrees are in $\Ip(n')$ for all $n'<n$, but $A$ might not be in $\Ip(n)$.
    
    As a consequence of those definitions, the relations $\Ip(n) \subseteq \Ap(n) \subseteq \Ep(n)$ hold for any integer $n>0$.
    So let $\Fp(n) = \Ap(n) - \Ip(n)$ be the set of \emph{forbidden} patterns in $\Ap(n)$: not forbidden as P2 patches but as part of a $3$-internal-regular induced subtree.
    We show by induction on $n$ that for $1 \leq n \leq 8$ the set $\Ip(n)$ is precisely the set appearing in Figure \ref{fig:poset} (brown tiles only), and $\Ip(n) = \emptyset$ for $n \geq 9$.
    Let $I' \in \Ip(n)$, where $n \geq 1$.

    \textsc{Case} $n = 1$.
    Then $I'$ is either a kite or a dart.
    Both tiles appear in Row 1 of Figure \ref{fig:poset} and they both admit at least one neighborhood of $3$ tiles.

    \textsc{Case} $n = 2$.
    For $n=2$, $\Ap(2)$ has five elements: the four brown subtrees shown in Row 2 of Figure \ref{fig:poset}, with their possible neighborhoods showing that they are in $\Ip(2)$, plus the brown subtree appearing in Figure \ref{fig:forbidden-2-chain} which cannot be in any $3$-internal-regular induced subtree.
    Indeed for each configuration $C$ in Figure \ref{fig:forbidden-2-chain}, if the brown tiles have degree $3$ then there is a cycle in $C$. 
    Hence the sets $\Ip(2)$ and $\Fp(2)$ are as illustrated in Row 2 of Figure \ref{fig:poset} and Figure \ref{fig:forbidden-2-chain} respectively.
\begin{figure}[h]
    \centering
    \includegraphics[scale=0.4]{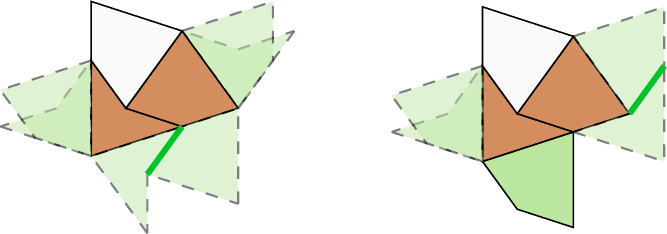}
    \caption{Up to isometry, the single element of $\Fp(2)$ (in brown) with its possible neighborhoods.}
    \label{fig:forbidden-2-chain}
\end{figure}

    \textsc{Case} $n = 3$.   For $n=3$,
     $\Ap(3)$ is the set of subtrees appearing in Row 3 of Figure \ref{fig:poset} or in Figure \ref{fig:forbidden-3-chains}.
    While all patterns in Figure \ref{fig:poset} are admissible, none of the five patches in Figure \ref{fig:forbidden-3-chains} is possible.
    Indeed in each patch of Figure \ref{fig:forbidden-3-chains}, at most one tile among each  pair of adjacent tiles $\{t_1, t_2\}$, $\{t_3, t_4\}$ and $\{t_5, t_6\}$ may belong to $I$, and by Lemma \ref{lem:n1n3} the tree should have 5 leaves.
    Let  $i_1$, $i_2$, $i_3$ be the chosen tiles in each of the three sets.
    Then there is a cycle formed by the tiles $i_1$, $i_2$, $i_3$, $t_7$ and $t_8$.
\begin{figure}[h]
    \centering
    \includegraphics[scale=0.4]{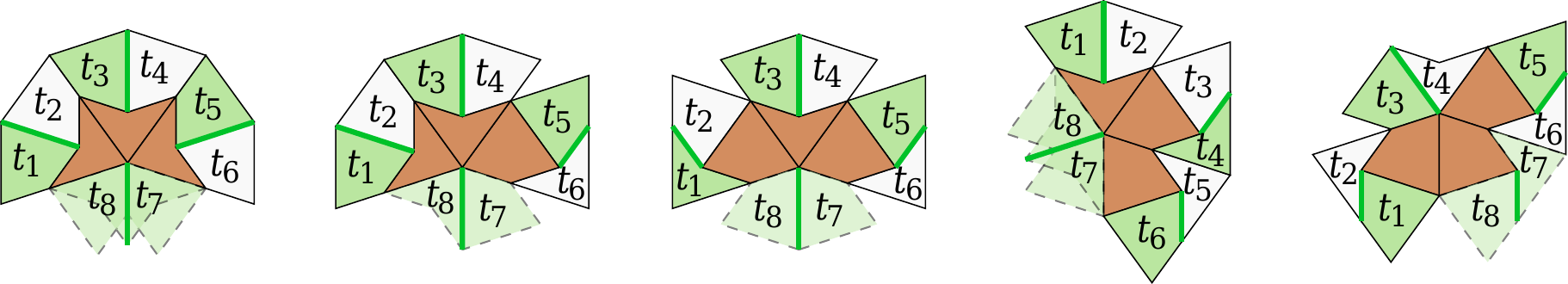}
\caption{Up to isometry, the five elements of $\Fp(3)$ (in brown), with their possible neighborhoods.}
    \label{fig:forbidden-3-chains}
\end{figure}

    \textsc{Case} $n = 4$. For $n=4$,
    observe that $I'$ cannot be a non-path subtree.
    Indeed, Figure \ref{fig:non-caterpillar-extensions-4} illustrates all non-path subtrees in $\Ep(4)$.
    Since each of them contains at least one forbidden pattern in $\Fp(3)$, we conclude that they cannot belong to $\Ap(4)$.
    Finally, all path subtrees in $\Ap(4)$ are depicted in Row 4 of Figure \ref{fig:poset}, so that $\Fp(4) = \emptyset$.
\begin{figure}[h]
    \centering
    \includegraphics[scale=0.4]{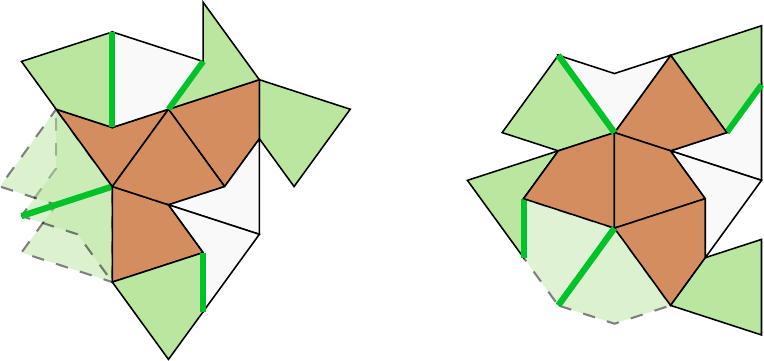}
\caption{Up to isometry, the non-path extensions in $\Ep(4)$.}
    \label{fig:non-caterpillar-extensions-4}
\end{figure}


    \textsc{Case} $n = 5$. For $n=5$,
     $\Ap(5)$ is the set of all isometries of the subtrees appearing in Figure \ref{fig:poset}, row $5$. Indeed, considering forced tiles around each element of $\Ip(4)$ and observations made before, no other extension is admissible.

    \textsc{Cases} $n = 6, 7, 8$.
    By inspecting all paths in $\Ep(n), n\in \{6,7,8\}$, one notices that $\Ap(n)$ is as depicted in Figure \ref{fig:poset}.
    Since they are all admissible, we conclude that $\Ap(n) = \Ip(n)$, so that $\Fp(n) = \emptyset$.

    \textsc{Case} $n = 9$.
    By inspecting all paths in $\Ep(9)$, we observe that none of them are in $\Ap(9)$, since they all contain some forbidden pattern in $\Fp(n')$ for some $n' < 9$.
    Hence, $\Ap(9) = \Ip(9) = \emptyset$.

    \textsc{Case} $n \geq 10$.
    By induction, since $\Ip(n - 1) = \emptyset$, $\Ep(n) = \emptyset$, hence $\Ap(n) = \emptyset$
    so that $\Ip(n) = \emptyset$, concluding the proof.
\end{proof}

Lemma \ref{lem:3-regular}  has an interesting consequence:

\begin{corollary}\label{corol:caterpillar}
Let $I$ be a $3$-internal-regular induced subtree of a P2-graph.
Then $I$ is a caterpillar.
\end{corollary}

\begin{proof}
    Let $I'$ be the derived tree of $I$.
    We need to show that $I'$ is a path graph.
    If $|I'| \leq 3$, then the result is immediate.
    Otherwise, assume that $|I'| \geq 4$ and, arguing by contradiction, that $I'$ is not a chain.
    Then $I'$ contains a subtree $J$ of $4$ tiles, such that one tile $t$ of $J$ verifies $\deg_{I'}(t) = 3$.
    But such a tile $t$ does not exist: in the proof of Lemma \ref{lem:3-regular}, we have shown that the set $\Ap(4)$ contains only path graphs.
    This contradicts our earlier assumption and concludes the proof.
\end{proof}

Before proving Proposition \ref{prop:borne-sup}, we introduce a useful definition.
\begin{definition}[Adapted from \cite{blondin-saturated}]
Let $I$ be an induced subtree of a P2 tiling and let $t_1, t_2$ be two adjacent tiles in $I$.
    We say that the pair of induced subtrees $(I_1,I_2)$ is a \emph{factorization of $I$ at $(t_1, t_2)$}, or that $I$ is the \emph{graft of $I_1$ and $I_2$} if $I = I_1 \cup I_2$, 
    $I_1 \cap I_2 = \{t_1, t_2\}$, $t_2$ is a leaf of $I_1$ and $t_1$ is a leaf of $I_2$.
    We then write $I = I_1 \comp_{(t_1,t_2)} I_2$ or simply $I = I_1 \comp I_2$ when $(t_1,t_2)$ is clearly identified.
\end{definition}
\begin{figure}[b]
    \centering
    \includegraphics[scale=2]{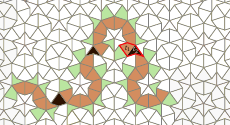}
    \caption{A factorization $I = I_1 \comp I_2$ at $(t_1,t_2)$. $I_1$ is the leftmost subtree of $I$, including $t_2$ as a leaf, and $I_2$ is the rightmost subtree of $I$, including $t_1$ as a leaf.}
    \label{fig:graft-example}
\end{figure}

The following identities are straightforward consequences of the previous definition:
For any induced subtrees $I_1$ and $I_2$ such that $I_1 \comp I_2$ exists,
\begin{eqnarray}
    n(I_1 \comp I_2) & = & n(I_1) + n(I_2) - 2, \\
    n_1(I_1 \comp I_2) & = & n_1(I_1) + n_1(I_2) - 2, \\
    n_2(I_1 \comp I_2) & = & n_2(I_1) + n_2(I_2), \\
    n_3(I_1 \comp I_2) & = & n_3(I_1) + n_3(I_2).
\end{eqnarray}

We are now ready to prove Proposition \ref{prop:borne-sup}.

\begin{proof}[of Proposition \ref{prop:borne-sup}]
    Exhaustive enumeration of all possible tile patterns shows that Inequation \ref{ineq:upper-bound} holds for $n \in \{0, 1, 2, 3\}$.
    Next, assume that $4 \leq n \leq 18$.
    On one hand, if $n$ is even, then by Lemma \ref{lem:loose-upper-bound}, we have $L_{P_{2}}(n) \leq L_{P_{2}}(2) + (n-2)/2 = 2 + (n-2)/2 = \lfloor n / 2 \rfloor + 1$.
    On the other hand, if $n$ is odd, then by Lemma \ref{lem:loose-upper-bound}, we have $L_{P_{2}}(n) \leq L_{P_{2}}(3) + (n-3)/2 = 2 + (n-1)/2 - 1$, hence $L_{P_{2}}(n) \leq \lfloor n / 2 \rfloor + 1$.

   We now consider the case $n \geq 19$.
    Arguing by contradiction, assume that there exists an integer $n\geq 19$ such that $L_{P_{2}}(n) > L_{P_{2}}(n - 17) + 8$.
    Choose $n$ minimal and let $I$ be a fully leafed induced subtree of size $n$.
  We claim that $I$ contains at least one tile of degree $2$.
    Indeed, if $n_2(I) = 0$, then
    $$19 \leq n(I)
         = n_1(I) + n_2(I) + n_3(I)
         = (n_3(I) + 2) + 0 + n_3(I)
         = 2n_3(I) + 2,$$
    which implies $n_3 \geq 9$, contradicting Lemma \ref{lem:3-regular}.
    Next, we observe that there exists a factorization $I = I_1 \comp I_2$ at $(t_1,t_2)$ such that $\deg_I(t_2) = 2$ and $t_2$ is the only tile of degree $2$ in $I_2$ (as in Figure \ref{fig:graft-example}).
    Therefore, it follows from Lemma \ref{lem:3-regular} that $n(I_2) \leq 19$: on the left side of $t_2$ there is only $t_1$ adjacent to $t_2$, and on its right there are at most 8 inner vertices (of degree 3) with as many leaves plus 1 at the end.
    To conclude, we consider two cases.

    First, assume that $n(I_2) \leq 18$.
    Since $n(I_2) = n_1(I_2) + n_2(I_2) + n_3(I_2)$, $n_2(I_2) = 1$ and $n_1(I_2) = n_3(I_2) + 2$, we have $n_1(I_2) = (n(I_2) + 1) / 2$, which also implies that $n(I_2)$ is odd.
    But by assumption $n_1(I) = L_{P_{2}}(n(I)) > L_{P_{2}}(n(I) - 17) + 8$ so
    $$n_1(I_1) = n_1(I) - n_1(I_2) + 2 > L_{P_{2}}(n(I) - 17) + 8 - n_1(I_2) + 2.$$ 
    Hence 
    $$n_1(I_1) > L_{P_{2}}(n(I) - 17) - (n(I_2) + 1) / 2 + 10 = L_{P_{2}}(n(I) - 17) + (19 - n(I_2)) / 2.$$
    Since $n(I_2)$ is odd, this rewrites as
    $$n_1(I_1) > L_{P_{2}}((n(I) - n(I_2) + 2) - (19 - n(I_2))) ~+~\lceil (19 - n(I_2)) / 2 \rceil $$
    Thus by Lemma \ref{lem:loose-upper-bound} we have $n_1(I_1) > L_{P_{2}}(n(I) - n(I_2) + 2)$
    i.e., $n_1(I_1) > L_{P_{2}}(n(I_1))$, which is impossible by definition of the leaf function.

    Finally, assume that $n(I_2) = 19$.
    Then $n_1(I_2) = 10$, $n_2(I_2) = 1$, $n_3(I_2) = 8$ and
    $$n(I_1) = n(I) - n(I_2) + 2 = n(I) - 19 + 2 = n(I) - 17.$$
    Also, we have 
    $$n_1(I_1) =  n_1(I) - n_1(I_2) + 2 = n_1(I) - 10 + 2 = n_1(I) - 8$$
    but then the assumption yields the same contradiction:
    $$n_1(I_1) > L_{P_{2}}(n(I) - 17) + 8 - 8 = L_{P_{2}}(n(I_1)).$$
\end{proof}

\begin{figure}[h]
    \centering
    \includegraphics[scale=2]{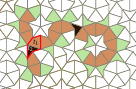}
    \caption{A fully leafed induced subtree which is not a caterpillar, but a graft of two caterpillars sharing tiles $t_1$ and $t_2$.}
    \label{fig:non-cat-flis-example}
\end{figure}

Before ending this section, it is worth mentioning that there exist fully leafed induced subtrees of tilings of type $P2$ that are not caterpillars. They are obtained by grafting a new caterpillar $C_{1}$ on a caterpillar $C$ with $8k$ tiles of degree $3$ and $k-1$ tiles of degree $2$. Indeed in order to increase the number of cells of $C$, a tile of degree $2$ must be added and this tile does not need to be at an extremity of $C'$ as illustrated in  Figure \ref{fig:non-cat-flis-example}.


\section{Caterpillars realizing the upper bound}

In this section we exhibit a family $(C_n)_{n\in\N}$ of fully leafed P2 caterpillars, that is caterpillars such that $n_1(C_n)=L_{P2}(n)$ for all $n\in\N$, thus proving:
\begin{proposition}\label{prop:optimal-cat}
Let $C_n$ be a Penrose caterpillar constructed as described in this section. 
Then
$$n_1(C_n) = 
\begin{cases}
    0 & \text{if } ~0\leq n\leq 1,\\
    \left\lfloor n/2 \right\rfloor +1 & \text{if } ~2\leq n\leq 19,\\
    n_1(C_{n-17})+8 & \text{if } ~n\geq 20.
\end{cases}
$$
\end{proposition}
The main difficulty is showing that these caterpillars are all allowed patches in P2 tilings.
For this purpose, we describe a recursive construction of optimal caterpillars, which can grow arbitrarily large.
We proceed by constructing local mappings applied on base caterpillars.

\subsection{The base caterpillars $\boldsymbol{C_{14}}$ and $\boldsymbol{C_{116}}$}

In addition to proving the upper bound, the previous section gives us optimal patches.
In particular, the 3-internal-regular caterpillars of internal size 8 (row $8$ of the graded poset in Figure \ref{fig:poset}) all include the same 3-internal-regular caterpillar of internal size 6 denoted $C_{14}$, whose derived path $C_{14}'$ is composed of $6$ kites surrounding a Star\footnote{The uppercase S is important, to avoid confusion with other stars which will appear later.}.
We use the name \emph{big sun} to refer to the patch consisting of a Star surrounded with ten kites (Figure  \ref{fig:kingdoms-P2}), since it is a \emph{sun} inflated once: each kite inflates into an \emph{ace}.
Since we are interested in leaves, i.e., in the tiles adjacent to the inner path of optimal caterpillars, we call a patch consisting of a big sun together with its adjacent tiles a \emph{flower} (Figure \ref{fig:7fleurs++}).

\begin{figure}[t]
    \centering
    \includegraphics[width=0.8\textwidth]{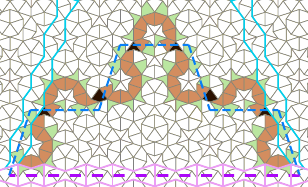}
    \caption{Optimal caterpillar $C_{116}$, consisting of seven 3-regular caterpillars of length 8, alternating with six isolated vertices of degree 2 (black), from which we removed the leaves adjacent to the endpoints of the inner path. Some leaves may be placed differently (as in Figure \ref{fig:poset}).}
    \label{fig:7fleurs}
\end{figure}

In order to have as many vertices of degree 3 as possible, each having an adjacent leaf as opposed to vertices of degree 2, an induced subtree should include as many copies  of flowers as possible. 
Ideally, we would like a caterpillar whose internal vertices (following the derived path) alternate 8 vertices of degree 3, then 1 vertex of degree 2, then 8 vertices of degree 3 and so on.
The question is, is it possible to find arbitrarily large such caterpillars in a $P2$ tiling?
Since Penrose tilings have the local isomorphism property, if such an arbitrarily large patch appears in one of them then it appears in all of them (infinitely many times).
Figure \ref{fig:7fleurs} shows that it is possible up to an induced subtree with 116 vertices, thus achieving the maximum of $L_{P2}(116)=56$ leaves.\\

Actually, we found optimal caterpillars much larger than this one but this specific subtree, which we denote $C_{116}$, is the base element of the following construction.
$C_{116}$ is  made of seven copies of $C_{14}$ (the shapes of the leaves being irrelevant) joined together by suitably adding 3 tiles: a black tile connects one ``extremal'' leaf of a copy to an ``extremal'' leaf of another copy, then a new leaf is added on each of them -- which thus becomes a vertex of degree 3 in the caterpillar.
It is placed ``on'' a Conway worm as in Figure \ref{fig:7fleurs}, consisting of long and short \emph{bowties} (the small patches delimited with pink segments) placed edge to edge.
In  Figure \ref{fig:7fleurs} and the following figures, this Conway worm is identified by a dashed purple line.
It is actually obtained when one inflates three times the tiles intersected by a dashed blue segment (then keeping only the tiles intersecting the purple line).

In fact, $C_{116}$ can be seen as a $\varphi^3$-decomposition of $C_{14}$.
Indeed, the three central copies of $C_{14}$ have their inner path in a \emph{king}'s kingdom, and most of their leaves in the worms surrounding it (directed by dashed blue segments).
When the king's kingdom and adjacent worms are inflated three times, we obtain the bottom (pink) worm as well as the partial worms going up on the sides (directed by the first and last blue segments) and all the tiles delimited by them -- plus others above, which are not needed here.
Hence whenever there is a copy of the top flower with the \emph{king} below the blue segment and worms on the sides in the same orientation, we can draw a copy of $C_{116}$ in the patch obtained with a $\varphi^3$-decomposition.
At each end of the caterpillar, it is always possible to extend the path on the side (through the kite-leaf) or below, thus crossing the purple segment (through the dart-leaf).

\subsection{Relation to Star tilings}

\begin{figure}
    \centering
    \includegraphics[width=0.8\textwidth]{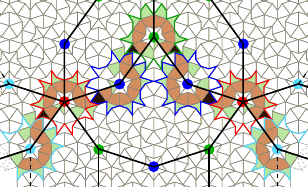}
    \caption{Base caterpillar $C_{116}$, with additional observations about the structure. Flowers are delimited in red, green, blue and cyan, their shapes determining the color of the central Star, i.e., of the vertices of the HBS shapes obtained by joining the centers of adjacent flowers.}
    \label{fig:7fleurs++}
\end{figure}

\begin{figure}[h]
    \begin{subfigure}[b]{0.35\textwidth}
        \centering
        \includegraphics[width=0.9\textwidth]{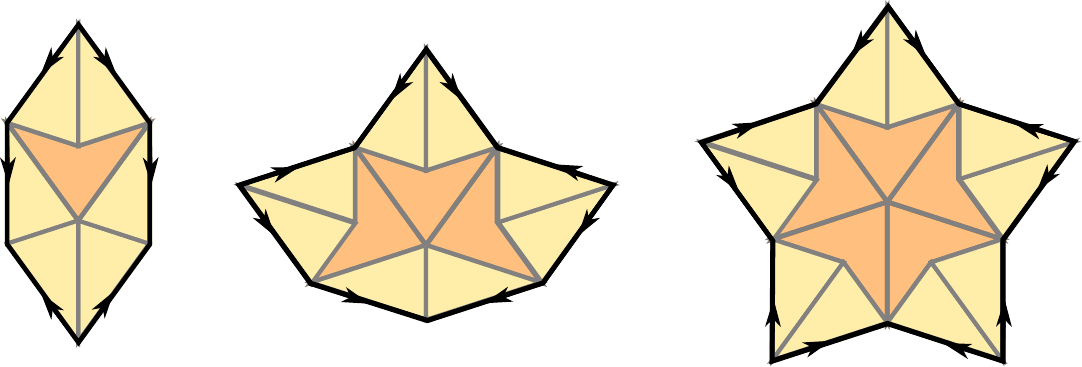}
        \caption{HBS tiles: hexagon, boat, star}
    \end{subfigure}
    \hfill
    \begin{subfigure}[b]{0.6\textwidth}
        \centering
        \includegraphics[width=\textwidth]{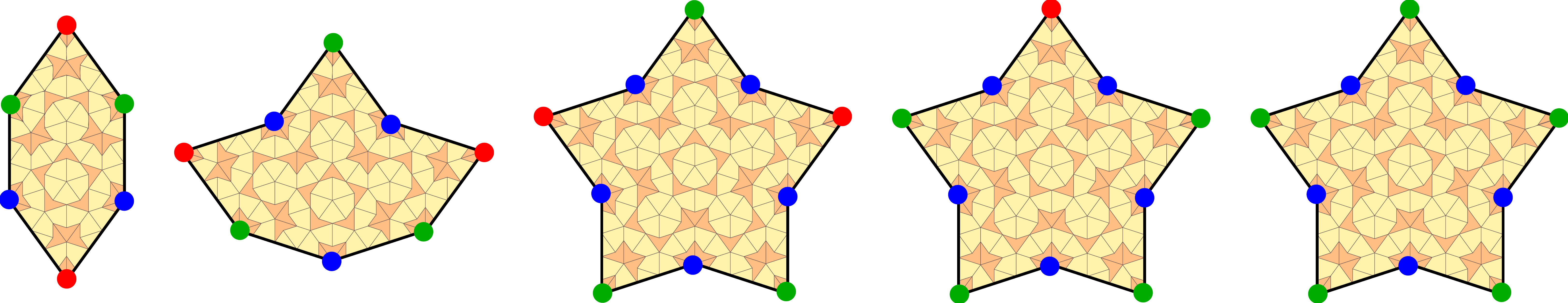}
        \caption{Star tiles: hexagon, boat, 2-star, 1-star and 0-star}
    \end{subfigure}
    \caption{HBS and Star tilesets, both decorated with kites and darts: the shapes are the same but the kites and darts which decorate the HBS tiles correspond to the $\varphi^3$-composition of those in the Star tileset.}
    \label{fig:star-tileset}
\end{figure}

Now let us add some information on the picture (Figure \ref{fig:7fleurs++}).
There are three types of Sun when we take into account their adjacent tiles \citep{bruijn1981}, yielding three types of flowers (after decomposition).
We attribute a color and an integer in $\{0,1,2 \}$ to each Star according to the fact that it is adjacent to either $0$ (red), $1$ (green) or $2$ (blue) Suns.
Whenever the optimal base patch $C_{116}$ appears, the ``inner'' flowers are always of the same type (02120, i.e., red-blue-green-blue-red), but the first and last ones that are met along the path are labelled either $1$ (green) or $2$ (blue), hence we color these Stars in cyan in Figure \ref{fig:7fleurs++}.
Notice that flowers in $C_{116}$ are adjacent to one another, and their adjacent tiles lie inside a short bowtie.
Hence the minimum distance between (the centers of) two ``adjacent'' \emph{Stars} is 3 long edges + 2 short edges, i.e., $3+2\varphi=\varphi^3$.
Joining two adjacent Star centers with a segment, one can see HBS shapes appear: hexagons, boats and stars\footnote{We use lowercase s for the star tiles, as opposed to the Star vertex configuration in P2 tilings.} (Figure \ref{fig:star-tileset}).
\begin{figure}[ht!]
    \centering
    \includegraphics[angle=90,width=0.8\textwidth]{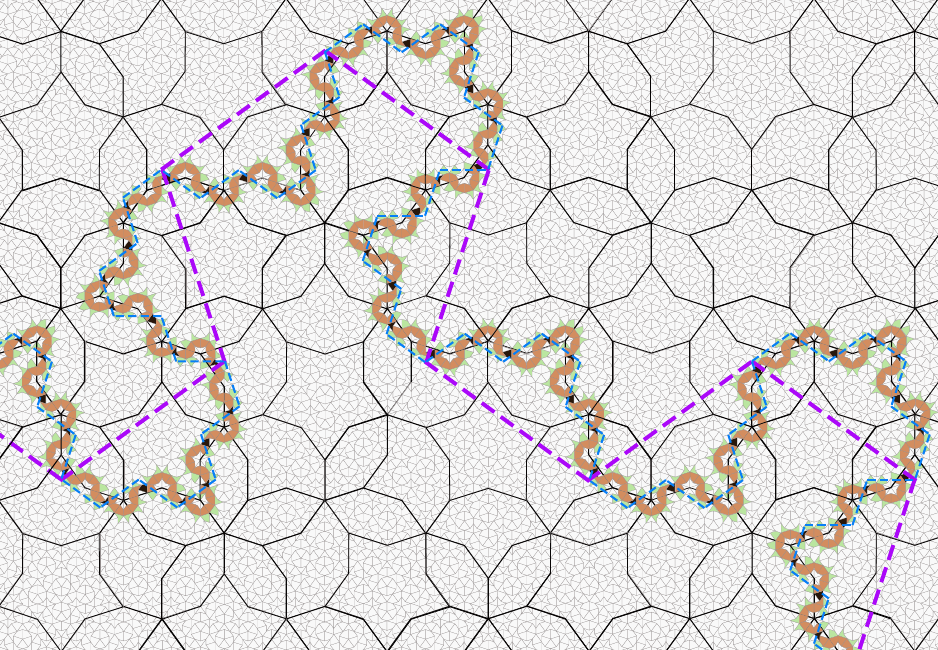}
    \caption{Long fully leafed caterpillar.}
    \label{fig:long-cat}
\end{figure}

HBS tilings are well known but, to our knowledge, this construction is new (\cite{Porrier_2024} presented a recent review of the existing literature).
The tiles obtained correspond to the $\varphi^3$-composition of the usual HBS tiles we would obtain by composing kites and darts.
In addition to this difference in decorations, by coloring the vertices of the HBS shapes according to the type of flower, the distribution of colors is unique on the vertices of hexagons and boats, while the stars come in three types, with either 0, 1 or 2 red vertices.
With this coloring of the vertices of the HBS shapes, we obtain what we call the \emph{Star tileset}, in which we call \emph{$i$-star} or \emph{star of type $i$} a star with $i$ red vertices (Figure \ref{fig:star-tileset}).
The sequence $02120$ followed by the derived path of $C_{116}$ appears only in boats (as shown in Figure \ref{fig:star-tileset}) and in stars of type 2.
As for the first and last flowers of $C_{116}$, due to the vertex configurations in Star tilings, 
they lie on the side of hexagons which share a 0 (red) vertex with either a boat or a star.

Figure \ref{fig:long-cat} shows how connected copies of $C_{116}$ meander around edges of a Star tiling.
The copies of $C_{116}$ in this long fully leafed caterpillar are connected in the same way as copies of  $C_{14}$ are connected to construct $C_{116}$.
One ``rule'' is forced: if a caterpillar is extended on the side, then the adjacent copy of $C_{116}$ is connected  from below (crossing its purple segment), and vice-versa.
This is due to the colors of the flowers that are connected, which have to be different (one blue and one green) to follow the side of a hexagon.
An interesting observation can be made: if we follow the line composed with dashed blue segments, the copies of $C_{14}$ alternate sides of the blue line: left-right-left-right-left and so on.
The same observation can be made about adjacent copies of $C_{116}$ as we follow the dashed purple line.
So there seems to be a logic of some sort.
Yet the path followed by the caterpillar still looks a bit random, with unexpected turns, so we use Star tilings to exhibit a larger structure.

\subsection{Optimal paths in Star tilings}\label{optimal-paths}

\begin{figure}[tb!]
    \centering
    \includegraphics[width=.9\textwidth]{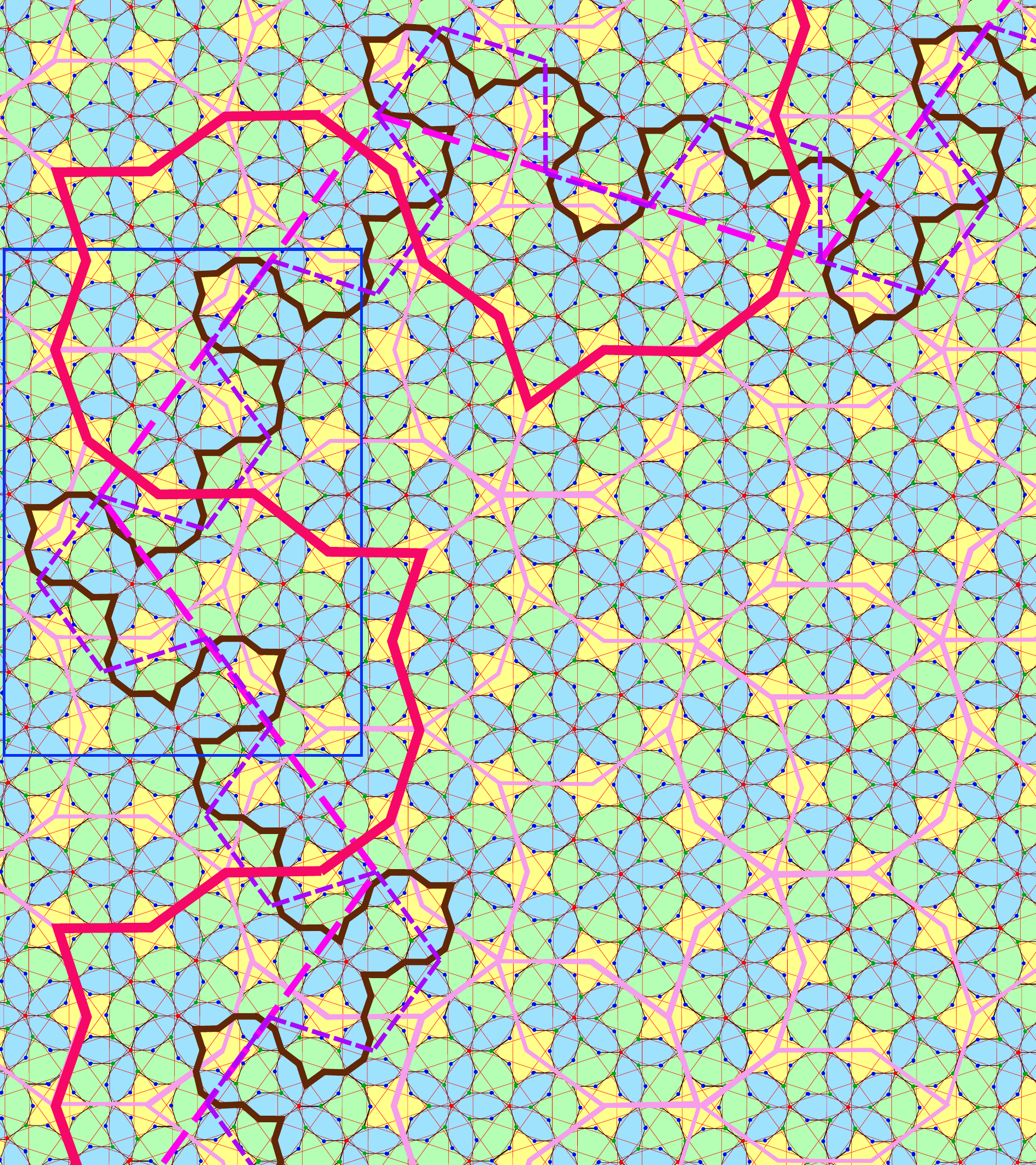}
    \caption{Long optimal path in Star tiling (in dark brown). The part of it which lies in the blue rectangle on the left border is the same as that of Figure \ref{fig:long-cat}, even though the underneath Star tiling is not the same. The red path is optimal in the $\varphi^3$-composition of the Star tiling (traced in pink).}
    \label{fig:long-path}
\end{figure}

In Figure \ref{fig:long-path}, the edges of the Star tiling around which the caterpillar meanders are traced in dark brown, but not the caterpillar itself: kites and darts are omitted since they are too small and we now know how to construct an optimal caterpillar from the appropriate path in the Star tiling.
We thus have an optimal path which is significantly longer than that of Figure \ref{fig:long-cat} -- and which includes it in the blue rectangle.
The same dashed purple segments are drawn, with endpoints always at the center of decagons (consisting of a boat and two hexagons) which are on both sides of the same boat-star (BS) configuration (a boat lying on the hollow part of a star).
Also, the $\varphi^3$-composition of the tiling appears in pink, along with an optimal path in it (in red).
It can be directly obtained by joining the centers of adjacent stars (instead of composing three times).

As dashed magenta segments are drawn beneath each $\varphi^3$-$C_{116}$ in the same way as the purple ones in the ``unit'' version $C_{116}$, 
it appears that each magenta segment can be decomposed into 7 purple segments forming a line from one endpoint of the magenta segment to the other, and the decomposition alternates sides around the magenta line, crossing the purple line at each turn.\\

\begin{figure}[bt]
    \centering
    \includegraphics[width=.85\textwidth]{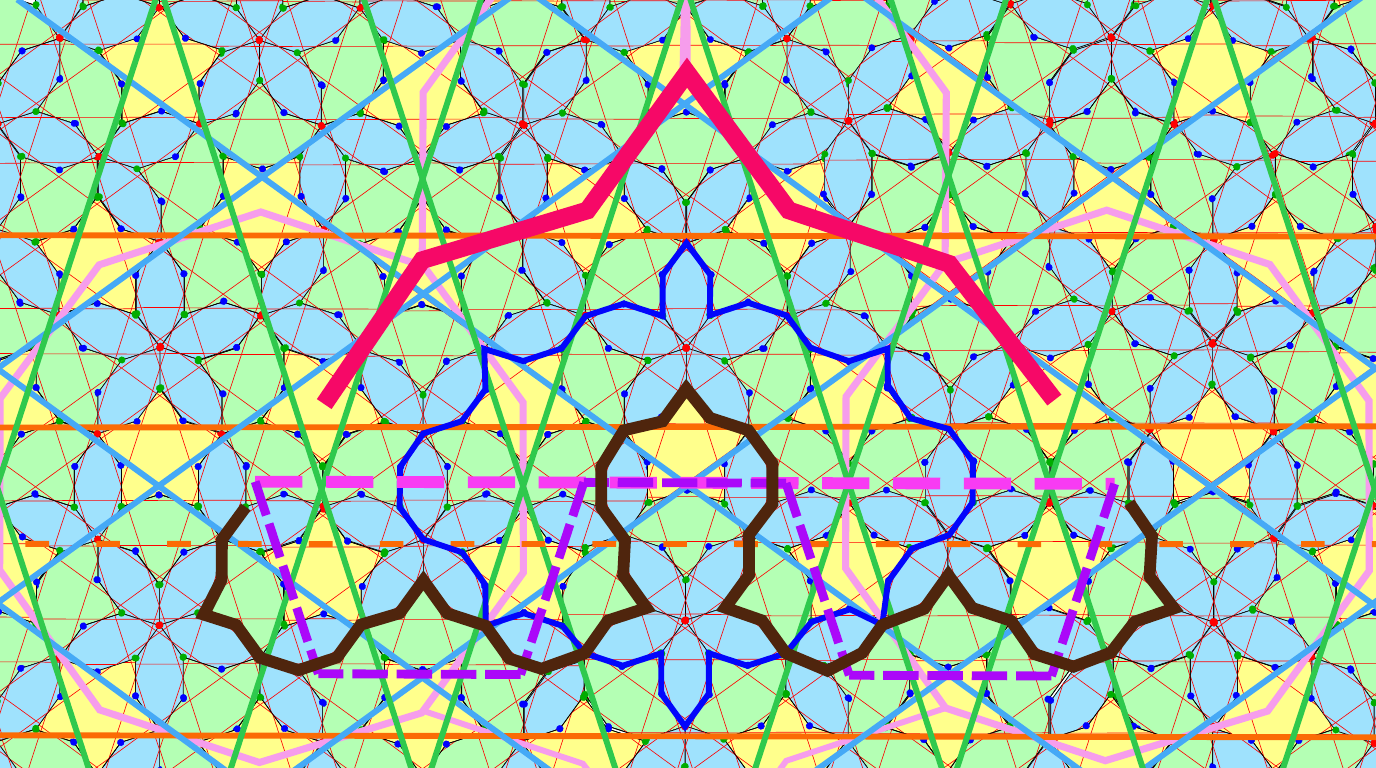}
    \caption{Decomposition of a single magenta segment and the corresponding (red) optimal path in the $\varphi^3$-Star tiling, yielding a $\varphi^3$-$C_{116}$. The 7 purple segments yield an optimal (dark brown) path in the Star tiling, which yields a fully leafed caterpillar $C_{830}$ -- composed with 7 copies of $C_{116}$ and 3 additional tiles at each of the 6 junctions (turns of the purple line).}
    \label{fig:path-decomp}
\end{figure}

In all dashed lines (blue, purple, magenta), the angles of the turns are $\pm 3\pi/5$ and the directions of the segments are those of the Ammann bars of the tiling (thin red lines).
The spacing between two consecutive parallel Ammann bars can take only two values: it is either long ($L$) or short ($S$).
The substitution rules are those of the Fibonacci word \citep{grunbaum2016}.
The blue, green and orange straight lines in Figure \ref{fig:path-decomp} are the $\varphi^3$-composition of the Ammann bars, with spacings denoted respectively $\mathcal{L}$ and $\mathcal{S}$.
Since the pattern is symmetrical through the central vertical line, we use the same color for angles $2\pi/5$ and $3\pi/5$ (green) as well as for $\pi/5$ and $4\pi/5$ (blue).

Figure \ref{fig:path-decomp} shows an isolated base portion of the path and its decomposition.
Actually, each dashed segment goes through the middle of a worm which can be flipped: the Ammann bar closest to the segment (in the same direction) could also be on the other side of the segment, at the same distance.
The impact of such a flip on the tiling is that the decagons and $BS$ configurations inside the Conway worm are reflected through the (dashed) central line.
The magic of our configurations is that the decomposition still works!
Indeed, when a $BS$ is in a Conway worm, the star is necessarily of type 2, thus there are hexagons on both sides of the $BS$, and when a worm is decomposed we still have a worm along the same dashed line in the decomposed tiling so the $BS$ configurations are the same. As a result, the middle of the decomposed path can follow the sides of the hexagons surrounding the $BS$, regardless of the orientation of the worm.

As for the part of the path below the worm, due to the presence of a 2-star and the tiles it forces around it (the kingdom of the star S2 is outlined in royal blue on Figure \ref{fig:path-decomp}), the configurations above and below the worm are the same, again whatever the orientation of the worm, and the path can therefore be the same.
Note that due to restrictions on Ammann bars, if the worm along the magenta segment is flipped then the one along the parallel purple segments cannot be flipped. As a result, the path is forced to be prolonged on the sides of the corresponding $C_{116}$'s and not below.
That is what is happening here with the $\varphi^3$-decomposition.



\subsection{Construction of the whole family $\boldsymbol{(C_n)_{n\in\N}}$}
\label{sectionCn}

We are now ready to describe a recursive procedure for defining a family of fully leafed caterpillars $(C_n)_{n\in\N}$.
We already have $C_{116}$ and its successive inflations, by the construction above.
The following algorithm allows us to construct caterpillars for all intermediate sizes (between inflations).
For a given $n\in\N$, we first find an inflation $C$ of $C_{116}$ with $m$ vertices where $m>n+15$.
Since $C$ is symmetric, we arbitrarily choose an endpoint of the derived path $C'$.
From there, we follow the path until we hit a vertex $v$ of degree 2.
This is the first of the new fully leafed induced subtree $T$ which the algorithm will output.
Following the path from $v$, we add one by one each vertex of the path and its adjacent leaf (if any) until $T$ has $n$ vertices.

\begin{algorithmic}[1]
  \Function{FLIS}{$n$: natural} : caterpillar
    \State $C \gets$ caterpillar of size $m$ obtained by inflations of $C_{116}$, such that $m>n+15$
    \State $v \gets$ first black vertex of $C$ (starting from an endpoint of $C'$)
    \State $T \gets$ empty caterpillar
    \State $k \gets 0$
    \While{$k<n$}
        \State \algorithmicif~$v$ is not yet in $T$ \algorithmicthen~add $v$ in $T$
        \State \textbf{else}\algorithmicif~$v$ has an adjacent leaf $\ell$ which is not yet in $T$ \algorithmicthen~add $\ell$ in $T$
        \State \textbf{else} $v \gets$ next vertex of $C'$ 
        \State $k \gets k+1$
    \EndWhile
    \State \Return $T$
  \EndFunction
\end{algorithmic}

\medskip
Since the derived path alternates 8 vertices of degree 3 (in the caterpillar) with 1 vertex of degree 2, and we take all the leaves, the number of leaves reaches the upper bound for each $n\in\N$.
We have proved Proposition \ref{prop:optimal-cat}. As a result, we have:
\begin{corollary}\label{prop:borne-inf}
The leaf function of the a P2 tiling satisfies the following inequality:
$$
L_{P2}(n) \geq 
\begin{cases}
    0 & \text{if } ~0\leq n\leq 1,\\
    \left\lfloor n/2 \right\rfloor +1 & \text{if } ~2\leq n\leq 18,\\
    L_{P2}(n-17)+8 & \text{if } ~n\geq 19.
\end{cases}
$$
\end{corollary}

Theorem \ref{thm1} follows from Proposition \ref{prop:borne-sup} and Corollary \ref{prop:borne-inf}.


\section{Final remarks and perspectives}

Though the recursive formula for $L_{P2}$  (Theorem \ref{thm1}) is quite simple, it can be more practical to have a closed formula, which we provide here.
\begin{proposition}\label{closed-formula}
The closed formula of the leaf function of P2 graphs is
$$L_{P2}(n) = 
\begin{cases}
    0 & \text{if } ~0\leq n\leq 1,\\
    \left\lfloor n/2\right\rfloor +1 & \text{if } ~2\leq n\leq 18,\\
    \left\lfloor (n+1)/2\right\rfloor  & \text{if } ~19\leq n\leq 35,\\
    8\left\lfloor n/17\right\rfloor + \left\lfloor (n\text{ \emph{mod} }17)/2\right\rfloor +\mathbb{1}(n\text{ \emph{mod} }17>1) & \text{if } ~n\geq 36.
\end{cases}
$$
\end{proposition}
\begin{proof}
We already know $L_{P2}(n)$ for $n\leq 18$.\\
For $19\leq n\leq 35$, we have $L_{P2}(n)=L_{P2}(n-17)+8$ with $2\leq n-17\leq 18$, hence
$$L_{P2}(n)=\left\lfloor (n-17)/2\right\rfloor+1+8 = \left\lfloor (n+1-18)/2\right\rfloor+9\left\lfloor (n+1)/2\right\rfloor.$$
We then show an intermediate result, by strong induction: for $n\geq36$,
$$L_{P2}(n)=L_{P2}\left(n-17\left(\left\lfloor n/17\right\rfloor-1\right)\right)+8\left(\left\lfloor n/17\right\rfloor-1\right).$$
\textsc{Base case}.  
    For $n=36$,
    $$L_{P2}\left(36-17\left(\left\lfloor 36/17\right\rfloor-1\right)\right)+8\left(\left\lfloor 36/17\right\rfloor-1\right) = L_{P2}\left(36-17\right)+8 = L_{P2}(36).$$
\textsc{Induction hypothesis}. For some integer $n\geq36$, for all integers $k$ such that $36\leq k\leq n$, 
    $$L_{P2}(k)=L_{P2}\left(k-17\left(\left\lfloor k/17\right\rfloor-1\right)\right)+8\left(\left\lfloor k/17\right\rfloor-1\right).$$
\textsc{Induction}. For $n+1$ we have
    $$L_{P2}(n+1) = L_{P2}(n+1-17)+8 = L_{P2}(n-16)+8$$ 
    with $n-16\geq19$ so we can use the induction hypothesis:
    \begin{align*}
        L_{P2}(n+1) &= L_{P2}\left(n-16-17\left(\left\lfloor \frac {n-16}{17}\right\rfloor-1\right)\right)+8\left(\left\lfloor \frac {n-16}{17}\right\rfloor-1\right) +8\\
                    &= L_{P2}\left(n-16-17\left(\left\lfloor \frac {n+1}{17}\right\rfloor-2\right)\right)+8\left(\left\lfloor \frac {n+1}{17}\right\rfloor-2\right) +8\\
                    &= L_{P2}\left(n+1-17\left(\left\lfloor \frac {n+1}{17}\right\rfloor-1\right)\right)+8\left(\left\lfloor \frac {n+1}{17}\right\rfloor-1\right).
    \end{align*}
Finally, to arrive at the result we note $n=17q+r$ with $q=\left\lfloor n/17\right\rfloor$ and $r=n\mod17$.
Thus we have
    \begin{align*}
        L_{P2}(n) &= L_{P2}\left(17q+r-17(q-1)\right)+8(q-1)\\
                  &= L_{P2}\left(r+17\right)+8(q-1)\\
                  &= L_{P2}\left(r\right)+8q
    \end{align*}
As a result, if $r=0$ or $r=1$ then $L_{P2}(n) = 8q= 8q+\left\lfloor r/2\right\rfloor +\mathbb{1}(r>1)$.\\
Otherwise $L_{P2}(n) = \left\lfloor r/2\right\rfloor +1+8q= 8q +\left\lfloor r/2\right\rfloor +\mathbb{1}(r>1)$.
\end{proof}
\medskip

\begin{corollary}
    The asymptotic growth of the leaf function for P2 graphs is $L_{P2}(n)\sim 8n/17$.
\end{corollary}
Since $8/17\simeq0,47$, we are roughly in the middle of the interval whose bounds were proved \citep{porrier2019}, i.e., $2\varphi n/(4\varphi + 1) \leq L_{P2}(n) \leq \lfloor n/2\rfloor+1 $, which implies $0,433 n \leq L_{P2}(n) \leq n/2+1 $.\\

Substitutive properties can be used to define arbitrarily large constructions in tilings.
In this way, we found an infinite family of fully leafed induced subtrees, from which we deduced the expression of the leaf function for P2-graphs.
Our construction relies on HBS shapes, which were widely used to model alloys and study their chemical properties (see \cite{Porrier_2024} for more details and references), relating it to quasicrystals.

One thing surprised us: we expected the leaf function $L_{P2}$ to be non-ultimately periodic and depending on the irrationality of $\varphi$, but it turned out not to be the case.
Note that the segments around which the caterpillars meander follow the Ammann bars.
We could also stick more closely to the segments, passing through the bottom of the boats instead of their top (to get from one ``endpoint'' of a boat to the other).
The drawback is that this is not possible when the $BS$ pattern is inverted, with the star ``on top'' and the boat upside down below.
However, apart from the leaves, which can be either a kite or a dart, and this alternative path along the boats, it seems that fully leafed Penrose subtrees are very often sub-caterpillars of the $(C_k)_{k\in\N}$ family described in Section \ref{sectionCn} or close to it.

We formalize these observations by stating two conjectures.
The first one is inspired by the article of \cite{blondin-saturated}, where the authors introduced the concept of \emph{saturated  polyform}, which can be naturally extended to the context of P2 tilings.
Indeed, since the leaf function $L_{P2}$ satisfies a linear recurrence, there exists a linear function $\overline{L_{P2}}(n)=8n/17+b$, where $b \in \R$, and a positive integer $N$ such that
$L_{P2}(n)\leq  \overline{L_{P2}}(n)$ for all $n\geq N$.
By choosing $b$ as small as possible, the linear function $\overline{L_{P2}}$ is unique. 
A Penrose tree $T$ is called \emph{saturated} if $n_{1}(T)=\overline{L_{P2}}(n(T))$. An \emph{appendix} is a fully leafed induced subtree $T$ of a $P2$ tiling such that $T'$ contains precisely one leaf of $T$-degree $2$ and at most two cells of  $T$-degree $3$.
We are now ready to state the first conjecture:

\begin{conjecture}\label{conj:caterpillar}
Let  $T$ be a fully leafed Penrose subtree. 
\begin{enumerate}[(i)]
\item If $T$ is saturated then $T $ is a caterpillar.
\item If $T$ is non saturated, then there exists a saturated caterpillar $T_1$ and an appendix $T_2$ such that $T=T_{1}\diamond T_{2}$.
\end{enumerate}
\end{conjecture} 

It is worth mentioning that, in the second case of Conjecture \ref{conj:caterpillar}, the resulting subtree $T$ can be either a non-saturated caterpillar (if the graft is localized on one of the extremity of $T_1$) or a non caterpillar (if the graft is localized on the interior of $T_{1}$, as in Figure \ref{fig:non-cat-flis-example}).

\begin{figure}[t]
    \centering
    \includegraphics[width=0.8\textwidth]{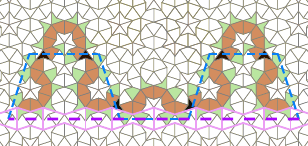}
    \caption{Optimal caterpillar, alternative to $C_{116}$, on the same worm.}
    \label{fig:boat-bottom}
\end{figure}

Moreover, we believe that the shape of any fully-leafed Penrose tree is closely related to the shape of the caterpillars belonging to the family $(C_n)_{n\in\N}$ in the following sense.
Let $f_{boat}$ be the transformation on the set $C(P2)$ of Penrose caterpillars such that $f_{boat}(C)$ is the caterpillar obtained from $C$ by substituting each subpath of $C'$ along the  hull of a boat (Figure \ref{fig:boat-bottom}) by the subpath along the top of the boat like in $C_{116}$. 
\begin{conjecture}
 Let $C$ be a Penrose caterpillar. Then $C$ belongs to one of the following two sets. 
 \begin{enumerate}[(i)]
\item $C$ belongs to the  graded poset $S$. 
\item  There exists a caterpillar $C_k$ such that the derived path $(f_{boat}(C))'$ is included in the derived path $(f_{boat}(C_{k}))'$, up to isomorphism.
\end{enumerate}
\end{conjecture}

In addition to enumerating all the fully leafed induced subtrees of P2-graphs, we can of course study the FLIS problem and look for the leaf functions of other tilings.
As such, despite the close links between Penrose tilings of types P2 and P3, we note that the leaf function is different for $P3$ as early as $n=12$: $L_{P3}(12)=6<7=L_{P2}(12)$.
This does not, however, preclude a possible transformation between fully leafed induced subtrees of either type, which would be interesting to explore.
Of course, the leaf function could be investigated  for the whole of Penrose's MLD class (notably HBS), and more generally for substitutive tilings.
In this respect, the study of different words derived from the $P2$ caterpillars could give us leads, potentially allowing us to apply results already known on Sturmian or billiard words.
We have started to investigate known words on $\{0,1,2\}$ corresponding to Star types, and words on angles of turns, to find out whether Penrose caterpillars' language is novel.

\acknowledgements 
We would like to thank {\it Cherry Arbor Design} for producing wooden kites and darts pieces, that were  useful in our journey into the Penrose leaf function.
We also thank the referee for their relevant suggestions.


\bibliographystyle{abbrvnat}
\bibliography{tilings,graphs}

\begin{thebibliography}{29}
\providecommand{\natexlab}[1]{#1}
\providecommand{\url}[1]{\texttt{#1}}
\expandafter\ifx\csname urlstyle\endcsname\relax
  \providecommand{\doi}[1]{doi: #1}\else
  \providecommand{\doi}{doi: \begingroup \urlstyle{rm}\Url}\fi

\bibitem[Abbott and Katchalski(1988)]{snake-in-the-box}
H.~Abbott and M.~Katchalski.
\newblock On the snake in the box problem.
\newblock \emph{Journal of Combinatorial Theory, Series B}, 45\penalty0
  (1):\penalty0 13--24, 1988.

\bibitem[Baake(1999)]{baake1999guide}
M.~Baake.
\newblock A guide to mathematical quasicrystals, 1999.
\newblock URL \url{https://arxiv.org/abs/math-ph/9901014}.

\bibitem[Baake and Grimm(2013)]{baake_grimm_2013}
M.~Baake and U.~Grimm.
\newblock \emph{Aperiodic Order. Vol 1. A Mathematical Invitation}, volume 149
  of \emph{Encyclopedia of Mathematics and its Applications}.
\newblock Cambridge University Press, Cambridge, 9 2013.
\newblock URL \url{http://oro.open.ac.uk/38933/}.

\bibitem[Blondin~Mass{\'e} et~al.(2018{\natexlab{a}})Blondin~Mass{\'e},
  de~Carufel, Goupil, Lapointe, Nadeau, and Vandomme]{blondin-induced}
A.~Blondin~Mass{\'e}, J.~de~Carufel, A.~Goupil, M.~Lapointe, {\'E}.~Nadeau, and
  {\'E}.~Vandomme.
\newblock Fully leafed induced subtrees.
\newblock In C.~Iliopoulos, H.~W. Leong, and W.-K. Sung, editors,
  \emph{Combinatorial Algorithms}, pages 90--101, Cham, 2018{\natexlab{a}}.
  Springer International Publishing.
\newblock ISBN 978-3-319-94667-2.
\newblock \doi{10.1007/978-3-319-94667-2_8}.

\bibitem[Blondin~Mass{\'e} et~al.(2018{\natexlab{b}})Blondin~Mass{\'e},
  de~Carufel, Goupil, and Samson]{blondin-poly}
A.~Blondin~Mass{\'e}, J.~de~Carufel, A.~Goupil, and M.~Samson.
\newblock Fully leafed tree-like polyominoes and polycubes.
\newblock In L.~Brankovic, J.~Ryan, and W.~F. Smyth, editors,
  \emph{Combinatorial Algorithms}, pages 206--218, Cham, 2018{\natexlab{b}}.
  Springer International Publishing.
\newblock ISBN 978-3-319-78825-8.

\bibitem[{Blondin Massé} et~al.(2018){Blondin Massé}, {de Carufel}, and
  Goupil]{blondin-saturated}
A.~{Blondin Massé}, J.~{de Carufel}, and A.~Goupil.
\newblock Saturated fully leafed tree-like polyforms and polycubes.
\newblock \emph{Journal of Discrete Algorithms}, 52-53:\penalty0 38--54, 2018.
\newblock ISSN 1570-8667.
\newblock \doi{https://doi.org/10.1016/j.jda.2018.11.004}.
\newblock URL
  \url{https://www.sciencedirect.com/science/article/pii/S1570866718301515}.
\newblock Combinatorial Algorithms – Special Issue Devoted to Life and Work
  of Mirka Miller.

\bibitem[Bodroza-Pantic et~al.(2013)Bodroza-Pantic, Pantic, Pantic, and
  Bodroza-Solarov]{bodroza2013}
O.~Bodroza-Pantic, B.~Pantic, I.~Pantic, and M.~Bodroza-Solarov.
\newblock Enumeration of hamiltonian cycles in some grid graphs.
\newblock \emph{MATCH Commun. Math. Comput. Chem.}, \penalty0 (70):\penalty0
  181--204, 2013.

\bibitem[Chan and Dill(1989)]{ChanDill1989}
H.~S. Chan and K.~A. Dill.
\newblock Compact polymers.
\newblock \emph{Macromolecules}, 22\penalty0 (12):\penalty0 4559--4573, 1989.
\newblock URL \url{https://doi.org/10.1021/ma00202a031}.

\bibitem[Chinnasamy et~al.(2019)Chinnasamy, Sivakumar, Selvakumari, and
  Suresh]{chinnasamy2019minimum}
A.~Chinnasamy, B.~Sivakumar, P.~Selvakumari, and A.~Suresh.
\newblock Minimum connected dominating set based rsu allocation for smartcloud
  vehicles in vanet.
\newblock \emph{Cluster Computing}, 22:\penalty0 12795--12804, 2019.

\bibitem[Collins et~al.(2017)Collins, Witte, Silverman, Green, and
  Gomes]{CWS2017}
L.~C. Collins, T.~G. Witte, R.~Silverman, D.~B. Green, and K.~K. Gomes.
\newblock Imaging quasiperiodic electronic states in a synthetic penrose
  tiling.
\newblock \emph{Nat Commun}, 8\penalty0 (15961), 2017.
\newblock URL \url{https://doi.org/10.1038/ncomms15961}.

\bibitem[{de Bruijn}(1981)]{bruijn1981}
N.~G. {de Bruijn}.
\newblock Algebraic theory of penrose's non-periodic tilings of the plane. i.
\newblock \emph{Indagationes Mathematicae (Proceedings)}, 84\penalty0
  (1):\penalty0 39--52, 1981.
\newblock ISSN 1385-7258.
\newblock \doi{https://doi.org/10.1016/1385-7258(81)90016-0}.
\newblock URL
  \url{https://www.sciencedirect.com/science/article/pii/1385725881900160}.

\bibitem[Diehl et~al.(2008)Diehl, Setyawan, and Curtarolo]{Diehl_2008}
R.~D. Diehl, W.~Setyawan, and S.~Curtarolo.
\newblock Gas adsorption on quasicrystalline surfaces.
\newblock \emph{Journal of Physics: Condensed Matter}, 20\penalty0
  (31):\penalty0 314007, jul 2008.
\newblock URL \url{https://dx.doi.org/10.1088/0953-8984/20/31/314007}.

\bibitem[Dubois(2000)]{DUBOIS2000}
J.-M. Dubois.
\newblock New prospects from potential applications of quasicrystalline
  materials.
\newblock \emph{Materials Science and Engineering: A}, 294-296:\penalty0 4--9,
  2000.
\newblock ISSN 0921-5093.
\newblock \doi{https://doi.org/10.1016/S0921-5093(00)01305-8}.
\newblock URL
  \url{https://www.sciencedirect.com/science/article/pii/S0921509300013058}.

\bibitem[Fernandes and Gouveia(1998)]{fernandes1998minimal}
L.~M. Fernandes and L.~Gouveia.
\newblock Minimal spanning trees with a constraint on the number of leaves.
\newblock \emph{European Journal of Operational Research}, 104\penalty0
  (1):\penalty0 250--261, 1998.

\bibitem[G\"ahler et~al.(1994)G\"ahler, Baake, and Schlottmann]{GBS1994}
F.~G\"ahler, M.~Baake, and M.~Schlottmann.
\newblock Binary tiling quasicrystals and matching rules.
\newblock \emph{Phys. Rev. B}, 50:\penalty0 12458--12467, Nov 1994.
\newblock URL \url{https://doi.org/10.1103/PhysRevB.50.12458}.

\bibitem[Gardner(1977)]{gardner1977}
M.~Gardner.
\newblock Mathematical games.
\newblock \emph{Scientific American}, 236\penalty0 (1):\penalty0 110--121,
  1977.
\newblock ISSN 00368733, 19467087.

\bibitem[Garel and Orland(1988)]{Garel_1988}
T.~Garel and H.~Orland.
\newblock Mean-field model for protein folding.
\newblock \emph{Europhysics Letters}, 6\penalty0 (4):\penalty0 307, jun 1988.
\newblock URL \url{https://dx.doi.org/10.1209/0295-5075/6/4/005}.

\bibitem[Grünbaum and Shephard(2016)]{grunbaum2016}
B.~Grünbaum and G.~C. Shephard.
\newblock \emph{Tilings and Patterns}.
\newblock Dover Publications, second edition, 2016.

\bibitem[Kweon et~al.(2022)Kweon, Kim, hwa Lee, Kim, and Lee]{KWEON2022117657}
J.~J. Kweon, H.-I. Kim, S.~hwa Lee, J.~Kim, and S.~K. Lee.
\newblock Quantitative probing of hydrogen environments in quasicrystals by
  high-resolution nmr spectroscopy.
\newblock \emph{Acta Materialia}, 226:\penalty0 117657, 2022.
\newblock ISSN 1359-6454.
\newblock \doi{https://doi.org/10.1016/j.actamat.2022.117657}.
\newblock URL
  \url{https://www.sciencedirect.com/science/article/pii/S1359645422000428}.

\bibitem[Li et~al.(1996)Li, Helling, Tang, and Wingreen]{LHTW1996}
H.~Li, R.~Helling, C.~Tang, and N.~Wingreen.
\newblock Emergence of preferred structures in a simple model of protein
  folding.
\newblock \emph{Science}, 273\penalty0 (5275):\penalty0 666--669, 1996.
\newblock \doi{10.1126/science.273.5275.666}.
\newblock URL
  \url{https://www.science.org/doi/abs/10.1126/science.273.5275.666}.

\bibitem[McGrath et~al.(2010)McGrath, Smerdon, Sharma, Theis, and
  Ledieu]{McGrath_2010}
R.~McGrath, J.~A. Smerdon, H.~R. Sharma, W.~Theis, and J.~Ledieu.
\newblock The surface science of quasicrystals.
\newblock \emph{Journal of Physics: Condensed Matter}, 22\penalty0
  (8):\penalty0 084022, feb 2010.
\newblock URL \url{https://dx.doi.org/10.1088/0953-8984/22/8/084022}.

\bibitem[Porrier(2024)]{Porrier_2024}
C.~Porrier.
\newblock {HBS Tilings Extended: State of the Art and Novel Observations}.
\newblock \emph{Electronic Proceedings in Theoretical Computer Science},
  403:\penalty0 156–163, June 2024.
\newblock ISSN 2075-2180.
\newblock URL \url{http://dx.doi.org/10.4204/EPTCS.403.32}.

\bibitem[Porrier and Blondin~Massé(2019)]{porrier2019}
C.~Porrier and A.~Blondin~Massé.
\newblock The leaf function for graphs associated with penrose tilings.
\newblock In \emph{2019 First International Conference on Graph Computing
  (GC)}, pages 37--44, 2019.
\newblock \doi{10.1109/GC46384.2019.00014}.

\bibitem[Porrier and Blondin~Massé(2020)]{PB2020}
C.~Porrier and A.~Blondin~Massé.
\newblock The {L}eaf function of graphs associated with {P}enrose tilings.
\newblock \emph{International Journal of Graph Computing}, 1.1:\penalty0 1--24,
  2020.
\newblock \doi{10.35708/GC1868-126721}.
\newblock URL
  \url{https://www.kspress.org/_files/ugd/e49175_b6e9d7522db94192a0e646ff5b5070ff.pdf}.

\bibitem[Robinson(2004)]{Robinson2004}
E.~A. Robinson.
\newblock Symbolic dynamics and tilings of $\mathbb{R}^{d}$.
\newblock \emph{Proceedings of Symposia in Applied Mathematics}, 60, 01 2004.
\newblock \doi{10.1090/psapm/060/2078847}.

\bibitem[Senechal(1995)]{Senechal1995}
M.~Senechal.
\newblock \emph{{Quasicrystals and Geometry}}.
\newblock Cambridge University Press, 1995.

\bibitem[Shioura et~al.(1997)Shioura, Tamura, and Uno]{SAT1}
A.~Shioura, A.~Tamura, and T.~Uno.
\newblock An optimal algorithm for scanning all spanning trees of undirected
  graphs.
\newblock \emph{SIAM Journal on Computing}, 26\penalty0 (3):\penalty0 678--692,
  1997.
\newblock \doi{10.1137/S0097539794270881}.
\newblock URL \url{https://doi.org/10.1137/S0097539794270881}.

\bibitem[Silva and Sivised(2018)]{desilva2018}
T.~N.~D. Silva and V.~Sivised.
\newblock A statistical mechanics perspective for protein folding from
  $q$-state potts model, 2018.
\newblock URL \url{https://arxiv.org/abs/1709.04813}.

\bibitem[Vedmedenko et~al.(2008)Vedmedenko, Mandel, and
  Lifshitz]{Vedmedenko01052008}
E.~Vedmedenko, S.~E.-D. Mandel, and R.~Lifshitz.
\newblock In search of multipolar order on the penrose tiling.
\newblock \emph{Philosophical Magazine}, 88\penalty0 (13-15):\penalty0
  2197--2207, 2008.
\newblock URL \url{https://doi.org/10.1080/14786430802302059}.

\end{thebibliography}

\end{document}